\documentclass[10pt]{article}
\usepackage{amssymb,amsmath,amsfonts,amsthm,enumerate}
\usepackage{subfigure}
\usepackage{graphicx}
\usepackage{epsfig}
\usepackage{float}
\usepackage[usenames]{color}

\def\no{\noindent} \def\p{\partial} \def\nb{\nonumber}
\def\Vh0{\stackrel{\circ}{V}_h} 
   
\def\Om{\Omega}   
\newcommand{\q}{\quad}
 
  \def\f{\frac}  
\def\D{\end{document}}   
   
\def\m{\mbox}   \def\Box{\sharp}

\newcommand{\lc}
{\mathrel{\raise2pt\hbox{${\mathop<\limits_{\raise1pt\hbox
{\mbox{$\sim$}}}}$}}}

\newcommand{\gc}
{\mathrel{\raise2pt\hbox{${\mathop>\limits_{\raise1pt\hbox{\mbox{$\sim$}}}}$}}}

\newcommand{\ec}
{\mathrel{\raise2pt\hbox{${\mathop=\limits_{\raise1pt\hbox{\mbox{$\sim$}}}}$}}}

\def\bb{\begin{equation}}  \def\ee{\end{equation}}

\def\beqn{\begin{eqnarray}}  \def\eqn{\end{eqnarray}}

\def\beqnx{\begin{eqnarray*}} \def\eqnx{\end{eqnarray*}}

\def\bn{\begin{enumerate}} \def\en{\end{enumerate}}

\def\bd{\begin{description}} \def\ed{\end{description}}
\def\bfg{\begin{figure}} \def\efg{\end{figure}}

\newtheorem{lemma}{Lemma}[section]
\newtheorem{example}{Example}[section]
\newtheorem{theorem}{Theorem}[section]

\newtheorem{algorithm}{Algorithm}[section]

 \def\x{{\bf x}}

\title{Quadratic Convergence of Levenberg-Marquardt Method \\
for Elliptic and Parabolic Inverse Robin Problems}
\author{
Daijun Jiang\footnote{School of Mathematics and Statistics $\&$ Hubei Key Laboratory of Mathematical
Sciences, Central  China Normal University, Wuhan, 430079, P.R.China. The
work of this author was financially supported by self-determined research funds
of CCNU from the colleges' basic research and operation of MOE (No. CCNU14A05039),
National Natural Science Foundation of China
(Nos. 11326233, 11401241 and 11571265). ({\tt jiangdaijun@mail.ccnu.edu.cn})}
\and
Hui Feng\footnote{School of Mathematics and Statistics, Wuhan University,
Wuhan 430072, China. The work of this author was supported by
National Natural Science Foundation of China (No.
91130022, No. 10971159 and No. 11161130003), The Doctoral
Fund of Ministry of Education of China (No. 20130141110026) and NCET
of China. ({\tt
hfeng.math@whu.edu.cn}). }
\and Jun
Zou\footnote{Department of Mathematics, The Chinese University of
Hong Kong, Shatin, Hong Kong. The work
of this author was substantially supported by Hong Kong RGC grants
(Projects 404611 and  405513). ({\tt zou@math.cuhk.edu.hk}).
}}

\begin{document}
\date{}
\maketitle

\begin{abstract}

We study the Levenberg-Marquardt (L-M) method for solving the highly nonlinear and ill-posed
inverse problem of identifying the Robin coefficients in elliptic and parabolic systems.
The L-M method transforms the Tikhonov regularized nonlinear non-convex minimizations
into convex minimizations. And the quadratic convergence of the L-M method
is rigorously established for the nonlinear elliptic and parabolic inverse problems for the first time,
under a simple novel adaptive strategy for selecting regularization parameters during the L-M iteration.
Then the surrogate functional approach is adopted
to solve the strongly ill-conditioned convex minimizations,
resulting in an explicit solution of the minimisation at each L-M iteration
for both the elliptic and parabolic cases.
Numerical experiments are provided to demonstrate the accuracy and efficiency of the methods.
\end{abstract}


%
\medskip
\section{Introduction}\label{sec:intro}
\setcounter{equation}{0}
We are concerned in this work with the determination of the Robin coefficient
in both stationary elliptic and time-dependent parabolic systems from
noisy measurement data on a partial boundary.
This is a highly nonlinear and ill-posed inverse problem and arises in many applications of practical importance.
The Robin coefficient may characterize the thermal properties of conductive materials
on the interface or certain physical processes near the boundary, e.g.,
it represents the corrosion damage profile
in corrosion detection \cite{fang04}\cite {ing97}, and indicates the
thermal property in quenching processes \cite {osm89}.

For the description of the model problems that are considered in this work,
we let $\Om\subset R^d~(d\geq1)$ be an open bounded and connected domain,
with a boundary $\p\Om$, which consists of two disjointed parts
$\p\Om=\Gamma_i\cup\Gamma_a$, with $\Gamma_i$ and $\Gamma_a$ each being
a $(d-1)$-dimensional polyhedral surface.
$\Gamma_i$ and $\Gamma_a$ are respectively the part
of the boundary that is inaccessible and accessible to experimental
measurements. Then we shall consider the inverse Robin problems associated with
the elliptic boundary value problem
\bb \left\{ \begin{array} {rclll}
-\nabla\cdot(a(\x)\nabla u)+c(\x) u &=&f(\x) & \m{in} &\Om\,, \\
a(\x)\frac{\p u}{\p n}+\gamma(\x) u&=&g(\x) & \m{on} &\Gamma_i\,, \\
a(\x)\frac{\p u}{\p n} &=&h(\x) & \m{on} &\Gamma_a\,,
\end{array}
\right. \label{ch1}
\ee
and the parabolic initial boundary value problem
\bb \left\{ \begin{array} {rclll}
\p_t u-\nabla\cdot(a(\x)\nabla u) &=&f(\x,t) & \m{in} &\Om\times [0,T]\,, \\
a(\x)\frac{\p u}{\p n}+\gamma(\x) u&=&g(\x,t) & \m{on} &\Gamma_i\times [0,T]\,, \\
a(\x)\frac{\p u}{\p n} &=&h(\x,t) & \m{on} &\Gamma_a\times [0,T]\,,\\
u(\x,0) &=&u_0(\x) & \m{in} &\Om\,.
\end{array}
\right. \label{ch1t}
\ee
The coefficients $a(\x)$ and $c(\x)$ are the heat conductivity and radiation coefficient,
satisfying that $\underline a\leq a(\x)\leq \bar a$  and
$\underline c\leq c(\x)\leq \bar c$  in $\Om$, where $\underline a,\, \bar a$ and
$\underline c,\, \bar c$ are positive constants. Functions $f$, $g$
and $h$ are the source strength,
ambient temperature and heat flux respectively. Both coefficients $\gamma(\x)$ in (\ref{ch1}) and (\ref{ch1t}) represent
the Robin coefficients, which will be the focus of our interest and is assumed to stay in the following feasible
constraint set:
$$
K:=\Big\{\gamma\in L^2(\Gamma_i);
\,\,0<\gamma_1\leq\gamma(\x)\leq\gamma_2 \,\,a.e.\,\, {\rm on}\,\, \Gamma_i\Big\},
$$
where $\gamma_1$ and $\gamma_2$ are two positive constants.
For convenience, we often write the solutions of the systems (\ref{ch1}) and (\ref{ch1t}) as $u(\gamma)$ to emphasize
their dependence on the Robin coefficient $\gamma$.

We are now ready to formulate the inverse problems of our interest in this work.

{\bf Elliptic Inverse Robin Problem}: recover the Robin coefficient $\gamma(\x)$ in (\ref{ch1}) on the inaccessible
part $\Gamma_i$  from the measurable data $z$ of $u$ on the accessible part $\Gamma_a$.

{\bf Parabolic inverse Robin problem}: recover the Robin coefficient $\gamma(\x)$ in (\ref{ch1t})
on the inaccessible part $\Gamma_i$ from the measurable data $z$ of $u$ on the accessible part $\Gamma_a$ over
the whole time range $[0,T]$.

The inverse Robin problems have been widely studied in literature; see \cite{fang04} \cite{cha04} \cite{jin07}\cite{jin12}\cite{jin11}
and the references therein.
The Gauss-Newton method  was applied in \cite{fang04} to solve the least-squares formulation of the
elliptic inverse Robin problem, but with no consideration of regularizations.
An $L^1$-tracking functional approach was suggested for the elliptic inverse Robin problem in \cite{cha04}.
Effectiveness and justifications of least-squares formulations with
regularizations were analysed in \cite{jin07}\cite{jin12}\cite{jin11} for the  Robin inverse problems,
and some iterative methods were applied to solve the resulting nonlinear least-squares minimizations.
However, we may observe a common feature of these existing methods,
which solve directly the nonlinear optimizations resulting from
least-squares formulations with regularisations, but these optimisation problems
are highly non-convex as the forward solution $u(\gamma)$ is nonlinear with respect to $\gamma$, and
strongly unstable at discrete level with fine mesh sizes and time step sizes due to the severe ill-posedness
of the inverse problems and the fact that noise is always present in the observation data.

In order to alleviate the effects of these drawbacks, we shall apply the L-M iterative method \cite{fan05} \cite{hanke97} \cite{lev44} \cite{mar63} \cite{yam01}
to solve the nonlinear optimizations resulting from least-squares formulations with regularisations
for the concerned inverse Robin problems. With the L-M method,
we need only to solve a convex optimization at each iteration.
Furthermore, in combination with the surrogate functional technique, we will not require the solution of
any optimisation problems in each iteration as the minimisers can be computed explicitly.
Another important novelty of this work is its establishment of the quadratic rate of convergence of the L-M method
for both the elliptic and parabolic inverse Robin problems. This appears to be the first time
in literature to demonstrate the quadratic convergence of the L-M method for a highly nonlinear ill-posed inverse problem.
Compared with general optimal control problems or general
direct nonlinear optimisation systems, the analysis on the quadratic rate of convergence of the L-M method here
is much more delicate and tricky, due to the severe ill-posedness, high nonlinearity and strong instability of the
current inverse problems and the direct effect on the convergence from two crucial parameters involved,
namely the regularization parameter and the noise level in the data.

The rest of the paper is organized as follows. In Section\,\ref{sec:uni},
we discuss the uniqueness of the nonlinear elliptic and parabolic inverse Robin problems.
In Sections\,\ref{sec:formulation} and \ref{sec:formulationt},  we formulate the Tikhonov
regularizations for the nonlinear elliptic and parabolic inverse Robin problem respectively and
study some mathematical properties of the resulting nonlinear optimisations.
In Subsections \ref{sub:b1} and \ref{sub:b2}, Fr\'{e}chet derivatives of the forward solution
of (\ref{ch1}) and (\ref{ch1t}) and corresponding adjoint operators are derived respectively.
In Subsections \ref{sec:linear} and \ref{sec:linearp}, the L-M iterative methods are formulated
and their quadratic convergences are established. The surrogate functional approach is
applied in Subsections\,\ref{sec:ddmrobin} and \ref{sec:ddmrobint} to solve the convex minimization at each L-M iteration
for the nonlinear elliptic and parabolic inverse Robin problem respectively.
Several numerical experiments are presented in Section\,\ref{sec:nsource} to illustrate the efficiency and
accuracy of the proposed methods. Some concluding remarks are given in
Section\,\ref{sec:conclu}.

Throughout this work, $C$ is often used for a generic positive constant. We shall use the symbol
$\langle \cdot, \cdot\rangle$ for the general inner product, and write the norms of the spaces
$H^m(\Omega)$, $L^2(\Omega)$,
$H^{{1}/{2}}(\Gamma)$ and $L^2(\Gamma)$ (for some $\Gamma\subset \p\Omega$)  respectively
as $\|\cdot\|_{m, \Omega}$, $\|\cdot\|_{\Omega}$, $\|\cdot\|_{1/2, \Gamma}$ and
$\|\cdot\|_{\Gamma}$.

\section{Uniqueness and local Lipschitz stability of the inverse Robin problems}\label{sec:uni}
\setcounter{equation}{0}
\setcounter{figure}{0}
\setcounter{table}{0}
In this section, we shall demonstrate the uniqueness and local Lipschitz stability of the Robin coefficients in the
concerned nonlinear elliptic and parabolic inverse Robin problems.

We first study the uniqueness and local Lipschitz stability of the elliptic inverse Robin problem.

\begin{theorem}\label{the:unique} {\bf (Uniqueness of the elliptic inverse Robin problem)}
Let $\gamma_1$ and $\gamma_2$ be two solutions to the elliptic inverse Robin problem as stated in Section\,\ref{sec:intro}, and
$\p\Om$ be $C^2$ smooth. Furthermore, we assume that ${\rm meas}(\{\x\in \Gamma_i; \,u(\x)=0\})=0$ where
$u$ is the solution to the forward system (\ref{ch1}), then
$\gamma_1=\gamma_2$ almost everywhere on $\Gamma_i$.
\end{theorem}
\no {\it Proof}. It is straightforward to verify using (\ref{ch1}) that $u(\gamma_1)-u(\gamma_2)$ satisfies
\beqn
\left\{ \begin{array} {rclll}
-\nabla\cdot(a(\x)\nabla (u(\gamma_1)-u(\gamma_2)))+
c(\x) (u(\gamma_1)-u(\gamma_2))&=&0 & \m{in} &\Om\,, \\
a(\x)\frac{\p (u(\gamma_1)-u(\gamma_2))}{\p n} &=&0 & \m{on} &\Gamma_a\,,\\
u(\gamma_1)-u(\gamma_2)&=&0 & \m{on} &\Gamma_a\,,
\end{array}
\right.
\label{ji1}
\eqn
and on the boundary $\Gamma_i$,
\bb
a(\x)\frac{\p (u(\gamma_1)-u(\gamma_2))}{\p n}+\gamma_1 u(\gamma_1)-\gamma_2 u(\gamma_2)=0.
\label{con1}
\ee
The unique continuation principle \cite{isa06} implies that $u(\gamma_1)-u(\gamma_2)=0$
in $\Om$. Hence, by the trace theorem and the weak form of the system (\ref{ji1}), we have
\beqnx
\|u(\gamma_1)-u(\gamma_2)\|_{\p\Om}\leq C\|u(\gamma_1)-u(\gamma_2)\|_{1,\Om}=0,
\eqnx
and for any $\varphi\in H^1(\Om)$,
\beqnx
\int_{\Gamma_i}a(\x)\frac{\p (u(\gamma_1)-u(\gamma_2))}{\p n}\varphi ds=\int_\Om a(\x)\nabla
(u(\gamma_1)-u(\gamma_2))\cdot\nabla \varphi + c(\x)(u(\gamma_1)-u(\gamma_2))\varphi d\x=0\,.
\eqnx
Therefore we immediately see sthat
\beqnx
a(\x)\frac{\p (u(\gamma_1)-u(\gamma_2))}{\p n}=0~~~~{\rm and}~~~~
 u(\gamma_1)= u(\gamma_2)~~~~{\rm on}~~~~\Gamma_i,
\eqnx
which, along with (\ref{con1}), leads to
 \beqnx
 u(\gamma_1)(\gamma_1-\gamma_2)=0~~~~{\rm on}~~~~\Gamma_i.
\eqnx
Now the assumption that ${\rm meas}(\{\x\in \Gamma_i:\,u(\x)=0\})=0$ implies
$\gamma_1=\gamma_2$ a.e. on $\Gamma_i$. \hfill $\Box$

For the local Lipschitz stability, let $\gamma^*$ be the true Robin coefficient,
we shall write for any positive constant $b$ that
\beqn
N(\gamma^*,b)=\{\gamma\in K;\,\,
\|\gamma-\gamma^*\|_{\Gamma_i}\leq b\}.
\label{cc1}
\eqn
Then we refer to \cite{jiang16} and give the following theorem to show the local Lipschitz stability.

\begin{theorem}\label{thm:locale}{\bf (Local Lipschitz stability of the elliptic inverse Robin problem)}
Assume
that $u(\gamma^*)\geq c_0>0$ on $\Gamma_i$, then there exists a positive constant $b$
such that the following stability estimate holds:
\bb
\|u(\gamma_1)-u(\gamma_2)\|_{\Gamma_a}\geq C\|\gamma_1-\gamma_2\|_{\Gamma_i} \q
\forall\, \gamma_1,\gamma_2\in N(\gamma^*,b)\,.
\label{cc5}
\ee
\end{theorem}

Next, we study the uniqueness and local Lipschitz stability of the parabolic inverse Robin problem.

\begin{theorem}\label{the:uniquep}{\bf (Uniqueness of the parabolic inverse Robin problems)}
Let $\gamma_1$ and $\gamma_2$ be two solutions to the parabolic inverse Robin problem
as stated in Section\,\ref{sec:intro}, and
$\p\Om$ be $C^2$ smooth. Moreover, we assume that ${\rm meas}(\{\x\in \Gamma_i; \,
u(\x, t)=0 ~\m{for} ~t\in (0,T)\})=0$, where $u$ is the solution to the forward system (\ref{ch1t}), then
$\gamma_1=\gamma_2$ almost everywhere on $\Gamma_i$.
\end{theorem}
\no {\it Proof}. It is straightforward to verify using (\ref{ch1t}) that $u(\gamma_1)-u(\gamma_2)$ satisfies
\beqn
\left\{ \begin{array} {rclll}
 \p_t(u(\gamma_1)-u(\gamma_2))-\nabla\cdot(a(\x)\nabla (u(\gamma_1)-u(\gamma_2)))&=&0 &
\m{in} &\Om\times [0,T]\,, \\
a(\x)\frac{\p (u(\gamma_1)-u(\gamma_2))}{\p n} &=&0 & \m{on} &\Gamma_a\times [0,T]\,,\\
u(\gamma_1)-u(\gamma_2)&=&0 & \m{on} &\Gamma_a\times [0,T]\,,\\
(u(\gamma_1)-u(\gamma_2))(\x,0)&=&0& \m{in} &\Om\,,
\end{array}
\right.
\label{ji2}
\eqn
and on the boundary $\Gamma_i\times [0,T]$,
\bb
a(\x)\frac{\p (u(\gamma_1)-u(\gamma_2))}{\p n}+\gamma_1 u(\gamma_1)-\gamma_2 u(\gamma_2)=0.
\label{pcon1}
\ee
The unique continuation principle \cite{isa06} implies that $u(\gamma_1)-u(\gamma_2)=0$
in $\Om\times [0,T]$. Hence,
by the trace theorem and the weak form of the system (\ref{ji2}), we have
\beqnx
\|u(\gamma_1)-u(\gamma_2)\|_{\p\Om}\leq C\|u(\gamma_1)-u(\gamma_2)\|_{1,\Om}=0,
\eqnx
and for any $\varphi\in L^2(0,T;H^1(\Om))$,
\beqnx
\int_0^T\int_{\Gamma_i}a(\x)\frac{\p (u(\gamma_1)-u(\gamma_2))}{\p n}\varphi ds dt&=&
\int_0^T\int_\Om a(\x)\nabla(u(\gamma_1)-u(\gamma_2))\cdot\nabla \varphi d\x dt\\
&& + \int_0^T\int_\Om\p_t(u(\gamma_1)-u(\gamma_2))\varphi d\x dt=0\,.
\eqnx
Therefore we immediately see that
\beqnx
a(\x)\frac{\p (u(\gamma_1)-u(\gamma_2))}{\p n}=0~~~~{\rm and}~~~~
 u(\gamma_1)= u(\gamma_2)~~~~{\rm on}~~~~\Gamma_i\times [0,T],
\eqnx
which, along with (\ref{pcon1}), yields sthat
 \beqnx
 u(\gamma_1)(\gamma_1-\gamma_2)=0~~~~{\rm on}~~~~\Gamma_i\times [0,T].
\eqnx
Now the assumption ${\rm meas}(\{\x\in \Gamma_i; \,u(\x, t)=0~{\rm for} ~t\in (0,T)\})=0$ implies
$\gamma_1=\gamma_2$ a.e. on $\Gamma_i$. \hfill $\Box$

For the local Lipschitz stability,  we also refer to \cite{jiang16} and give the following theorem
to show the local Lipschitz stability.
\begin{theorem}\label{thm:localp}{\bf (Local Lipschitz stability of the parabolic inverse Robin problem)}
Assume
that $u(\gamma^*)\geq c_0>0$ on $\Gamma_i\times [0,T]$, then there exists a constant $b>0$
such that
\bb
\|u(\gamma_1)-u(\gamma_2)\|_{L^2(0,T;L^2(\Gamma_a))}\geq C\|\gamma_1-\gamma_2\|_{\Gamma_i}\,
\quad \forall\, \gamma_1,\gamma_2\in N(\gamma^*,b)\,.
\label{ccc6}
\ee
\end{theorem}

\section{Elliptic inverse Robin problem and its L-M solution}\label{sec:formulation}
\setcounter{equation}{0}
\setcounter{figure}{0}
\setcounter{table}{0}

\subsection{Tikhonov regularization for elliptic inverse Robin problem}
\label{sub:b1}
In this section we first formulate the Levenberg-Marquardt method for solving
the nonlinear non-convex optimisation problems resulting from the least-squares formulation
of the elliptic inverse Robin problem as stated in Section\,\ref{sec:intro}, incorporated with Tikhonov regularization
to handle its ill-posedness and instability due to the presence of the noise in the observation data  \cite{jin11}.
We assume the noise level in the observation data $z^\delta$ of the true solution $u$
to the elliptic system (\ref{ch1})
is of order $\delta$, namely
\bb
\|u(\gamma^*)-z^\delta\|_{\Gamma_a}\leq \delta\,,
\label{ch2}
\ee
where $\gamma^*$ is the true Robin coefficient.
The elliptic inverse Robin problem is frequently transformed
into the following stabilized minimization system with Tikhonov regularization:
\beqn
\min_{\gamma\in K}{\cal J}(\gamma)
&=& \|u(\gamma)-z^\delta\|_{\Gamma_a}^2+\beta\|\gamma\|_{\Gamma_i}^2\,,
\label{ch3}
\eqn
where $\beta$ is the regularization parameter.
The formulation (\ref{ch3}) was shown to be stable in the sense that its minimizer depends continuously on the change of
the noise in the data $z^\delta$ \cite{jin11}.

%
For the subsequent analysis on the convergence of the Levenberg-Marquardt method
for solving the optimisation (\ref{ch3}), we shall frequently need the Fr\'{e}chet derivative of the forward
solution $u(\gamma)$ of system (\ref{ch1}). Let $w:=u'(\gamma)d$ be the Fr\'{e}chet derivative
at direction $d$, then it solves the following system:
\bb \left\{ \begin{array} {rccll}
-\nabla\cdot(a(\x)\nabla w)+c(\x) w &=&0 & \m{in} &\Om\,, \\
a(\x)\frac{\p w}{\p n}+\gamma \,w&=&-d \,u(\gamma) & \m{on} &\Gamma_i\,, \\
a(\x)\frac{\p w}{\p n} &=&0 & \m{on} &\Gamma_a\,.
\end{array}
\right. \label{ch6}
\ee
Let $u'(\gamma)^*$ be the adjoint operator of the Fr\'{e}chet derivative $u'(\gamma)$, then it it easy to verify that
$w^*:=u'(\gamma)^*p\in H^1(\Om)$ at a general direction $p$
solves the following system
\bb \left\{ \begin{array} {rccll}
-\nabla\cdot(a(\x)\nabla w^*)+c(\x) w^* &=&0 & \m{in} &\Om\,, \\
a(\x)\frac{\p w^*}{\p n}+\gamma \,w^*&=& 0& \m{on} &\Gamma_i\,, \\
a(\x)\frac{\p w^*}{\p n} &=&-p \,u(\gamma) & \m{on} &\Gamma_a\,.
\end{array}
\right. \label{ch7}
\ee
The following lemma gives an important relation for our later study.
\begin{lemma}\label{lem:adjoint}
The following relation holds for any directions $d$ and $p$:
\bb
\langle
w,\,u(\gamma)p\rangle_{\Gamma_a}
=\langle u(\gamma)d,\,w^*\rangle_{\Gamma_i}.
\label{ch8}
\ee
\end{lemma}
\no {\it Proof}. For any $\varphi,\psi\in H^1(\Om)$, we can readily derive the variational
forms of systems (\ref{ch6}) and (\ref{ch7}):
\bb
\int_\Om a(\x)\nabla
w\cdot\nabla \varphi d\x+\int_\Om c(\x)w\varphi d\x=
\int_{\Gamma_i}(-du(\gamma)-\gamma w)\varphi ds,
\label{ch9}
\ee
\bb
\int_\Om a(\x)\nabla
w^*\cdot\nabla \psi d\x+\int_\Om c(\x)w^*\psi d\x=
-\int_{\Gamma_i}\gamma w^*\psi ds-\int_{\Gamma_a}p u(\gamma)\psi ds.
\label{ch10}
\ee
Now (\ref{ch8}) follows by taking $\varphi=w^*$ and $\psi=w$ respectively in (\ref{ch9}) and (\ref{ch10}).
\hfill $\Box$

\subsection{Levenberg-Marquardt method and its convergence}\label{sec:linear}
The nonlinearity of the forward solution $u(\gamma)$ of the system (\ref{ch1})
makes the minimization (\ref{ch3}) highly nonlinear and non-convex
with respect to the Robin coefficient $\gamma$, as well as strongly unstable
at discrete level with fine mesh sizes and time step sizes due to the severe ill-posedness
of the inverse problem and the fact that noise is always present in the observation data.
To alleviate these difficulties in numerical solutions, we shall apply the Levenberg-Marquardt method
to solve (\ref{ch3}). For a given $\bar \gamma\in K$, we apply the linearization
\beqnx
u(\gamma)\approx u(\bar\gamma)+u'(\bar \gamma)(\gamma-\bar\gamma),
\eqnx
then we may solve the minimization system (\ref{ch3}) by the following Levenberg-Marquardt
iteration, which is widely used for general nonlinear optimization
problems \cite{fan05} \cite{yam01}:
\beqn
J(\gamma^{k+1})=\min_{\gamma\in K}J(\gamma)
:=\|u'(\gamma^k)(\gamma-\gamma^k)-(z^\delta-u(\gamma^k))\|_{\Gamma_a}^2
+\beta_k\|\gamma-\gamma^k\|_{\Gamma_i}^2\,.
\label{ch4}
\eqn

Before our study of the convergence of the iteration (\ref{ch4}),
we shall develop some auxiliary results.
\begin{lemma}\label{lem:estimation}
Assume the forward operator $u(\gamma)$ of system (\ref{ch1}) satisfies
that $u(\gamma)\in L^{\infty}(\Gamma_i)$ for $\gamma\in K$, then
there exist two positive constants $L$ and $c_1$ such that the following estimates hold
for all $\gamma,\,\bar\gamma\in K$:
\beqn
\|u(\gamma)-u(\bar \gamma)\|_{\Gamma_a}&\leq&L\|\gamma-\bar \gamma\|_{\Gamma_i},\label{es1}\\
\|u'(\bar \gamma)(\gamma-\bar \gamma)-(u(\gamma)-u(\bar \gamma))\|_{\Gamma_a}
&\leq&c_1\|\gamma-\bar \gamma\|^2_{\Gamma_i}\,.\label{es2}
\eqn
\end{lemma}
\no {\it Proof}. From the variational form of the system (\ref{ch1}), we can easily find that
\beqn
&&\int_\Om a(\x)\nabla (u(\gamma)-u(\bar \gamma))\cdot\nabla \varphi d\x+
\int_\Om c(\x)(u(\gamma)-u(\bar \gamma))\varphi d\x\nb\\
&+&\int_{\Gamma_i}\bar\gamma (u(\gamma)-u(\bar \gamma)) \varphi ds=
-\int_{\Gamma_i}(\gamma-\bar\gamma)u(\gamma)\varphi ds.
\label{es3}
\eqn
Taking $\varphi=u(\gamma)-u(\bar \gamma)$ and using the lower bounds
of $a(\x)$, $c(\x)$, and $\bar\gamma$, we derive
\beqnx
&&\min\{\underline a,\underline c\}\|u(\gamma)-u(\bar \gamma)\|^2_{1,\Om}+\gamma_1\|u(\gamma)-u(\bar \gamma)\|^2_{\Gamma_i}
\leq|\int_{\Gamma_i}(\gamma-\bar\gamma)u(\gamma)(u(\gamma)-u(\bar \gamma)) ds|\\
&\leq&\|u(\gamma)\|_{L^{\infty}(\Gamma_i)}\|\gamma-\bar\gamma\|_{\Gamma_i}
\|u(\gamma)-u(\bar \gamma)\|_{\Gamma_i}\leq C\|\gamma-\bar\gamma\|_{\Gamma_i}
\|u(\gamma)-u(\bar \gamma)\|_{\Gamma_i}.
\eqnx
Then it follows by the Cauchy-Schwarz inequality that
\beqnx
\min\{\underline a,\underline c\}\|u(\gamma)-u(\bar \gamma)\|^2_{1,\Om}\leq \frac{C^2}{2\gamma_1}\|\gamma-\bar\gamma\|^2_{\Gamma_i}\,.
\eqnx
Now estimate (\ref{es1}) follows directly from this inequality and the trace theorem.
To verify the estimate (\ref{es2}), we first show
\bb
\|u'(\gamma)d\|_{\p\Om}\leq C\|d\|_{\Gamma_i}.
\label{es15}
\ee
Indeed,  choosing
$\varphi=u'(\gamma)d$ in (\ref{ch9}), we readily get
\beqnx
\int_\Om a(\x)|\nabla
(u'(\gamma)d)|^2 d\x+\int_\Om c(\x)|u'(\gamma)d|^2 d\x
+\int_{\Gamma_i}\gamma |u'(\gamma)d|^2 ds=
-\int_{\Gamma_i}du(\gamma)(u'(\gamma)d) ds.
\eqnx
Then it follows by the Cauchy-Schwarz inequality that
\beqnx
\min\{\underline a,\underline c\}\|u'(\gamma)d\|^2_{1,\Om}\leq \frac{C^2}{2\gamma_1}\|d\|^2_{\Gamma_i},
\eqnx
which, along with the trace theorem, gives (\ref{es15}) immediately.

Next, we prove the estimate (\ref{es2}).
Taking $\gamma=\bar\gamma$ and $d=\gamma-\bar\gamma$ in (\ref{ch9}), we have
\beqn
&&\int_\Om a(\x)\nabla u'(\bar\gamma)(\gamma-\bar\gamma)\cdot\nabla \varphi d\x+
\int_\Om c(\x)u'(\bar\gamma)(\gamma-\bar\gamma)\varphi d\x\nb\\
&+&\int_{\Gamma_i}\bar\gamma
u'(\bar\gamma)(\gamma-\bar\gamma)\varphi ds=
-\int_{\Gamma_i}(\gamma-\bar\gamma)u(\bar\gamma)\varphi ds ~~\forall\,\varphi\in H^1(\Om)\,.
\label{es5}
\eqn
Subtracting (\ref{es5})  from (\ref{es3}) yields
\beqnx
&&\int_\Om a(\x)\nabla (u'(\bar \gamma)(\gamma-\bar \gamma)-(u(\gamma)-u(\bar \gamma)))\cdot\nabla \varphi d\x+
\int_{\Gamma_i}\bar\gamma (u'(\bar \gamma)(\gamma-\bar \gamma)-(u(\gamma)-u(\bar \gamma))) \varphi ds\nb\\
&+&\int_\Om c(\x)(u'(\bar \gamma)(\gamma-\bar \gamma)-(u(\gamma)-u(\bar \gamma)))\varphi d\x=
\int_{\Gamma_i}(\gamma-\bar\gamma)(u(\gamma)-u(\bar\gamma))\varphi ds.
\eqnx
Then applying the trace theorem, Lagrange mean value theorem and inequality (\ref{es15}), we derive
\beqnx
&&\|u'(\bar \gamma)(\gamma-\bar \gamma)-(u(\gamma)-u(\bar \gamma))\|_{\Gamma_a}
\leq C\|u'(\bar \gamma)(\gamma-\bar \gamma)-(u(\gamma)-u(\bar \gamma))\|_{1,\Om}\\
&\leq&C\|(\gamma-\bar\gamma)(u(\gamma)-u(\bar\gamma))\|_{\Gamma_i}=
C\|(\gamma-\bar\gamma)u'(\xi)(\gamma-\bar\gamma)\|_{\Gamma_i}\leq c_1\|\gamma-\bar\gamma\|^2_{\Gamma_i},
\eqnx
where $\xi$ is some element in $K$ between $\gamma$ and $\bar\gamma$.   $\Box$

Now we are ready to establish a quadratic rate on the convergence of the Levenberg-Marquardt method
(\ref{ch4}), under the following basic condition:
\bb
\|u(\gamma)-u(\gamma^*)\|_{\Gamma_a}\geq c_2\|\gamma-\gamma^*\|_{\Gamma_i} \q \forall\,\gamma\in N(\gamma^*,b)
\label{es8}
\ee
where $c_2$ and $b$ are two positive constants with $b\in (0,1)$.
Here $N(\gamma^*,b)$ denotes the $b$-neighborhood of the true Robin coefficient $\gamma^*\in K$ defined
in (\ref{cc1}).
Assumption (\ref{es8}) is the frequently adopted basic condition
to ensure the quadratic convergence of the Levenberg-Marquardt method
for most direct nonlinear optimization problems \cite{fan05}  \cite{yam01}, so it is natural to bring it
to the current nonlinear ill-posed inverse problems. The condition (\ref{es8}) may be viewed as a direct motivation
of the local Lipschitz stability of the elliptic inverse Robin problems (see estimate (\ref{cc5}) in Theorem \ref{thm:locale}).

It is a well-known technical difficulty in a practical numerical realisation of any Tikhonov regularised optimisation system like
the ones \eqref{ch3} and \eqref{ch4} to choose a reasonable and effective regularization parameter $\beta$ or $\beta_k$.
Another important novelty of this work is our suggestion of
a very simple and easy implementable choice of the parameter $\beta_k$ based on the following rule:
\bb
\beta_k=\|u(\gamma^k)-z^\delta\|^2_{\Gamma_a}\,.
\label{es7}
\ee
And surprisingly, as we shall demonstrate below, this choice of the regularization parameter $\beta_k$ ensures
a quadratical convergence of the resulting Levenberg-Marquardt iteration (\ref{ch4}).

Considering the presence of the noise (see \eqref{ch2}), it is reasonable for us to terminate the L-M iteration (\ref{ch4})
when its minimizer $\gamma^k$ is accurate enough in terms of the noise level, more specifically,
we shall terminate the iteration if the following criterion is realised:
\bb
c_2\|\gamma^k-\gamma^*\|_{\Gamma_i}<2\delta \q \m{or} \q
\|u(\gamma^k)-z^\delta\|_{\Gamma_a}<\sqrt{\delta}\,.
\label{aj1}
\ee

\begin{lemma}\label{lem:oestimation}
Under the conditions (\ref{es8}), \eqref{es7} and (\ref{aj1}), if $\gamma^k\in N(\gamma^*,b)$ then
$\gamma^{k+1}$ generated by the iteration (\ref{ch4}) satisfies
\beqn
\|u'(\gamma^k)(\gamma^{k+1}-\gamma^k)-(z^\delta-u(\gamma^k))\|_{\Gamma_a}&\leq&
c_3(\|\gamma^{k}-\gamma^*\|^2_{\Gamma_i}+\delta), \label{es11}\\
\|\gamma^{k+1}-\gamma^k\|^2_{\Gamma_i}&\leq& c_4(\|\gamma^{k}-\gamma^*\|^2_{\Gamma_i}+\delta),
\label{es18}
\eqn
where constants $c_3$ and $c_4$ are given explicitly by $c_3=\sqrt{\max\{2c_1^2+2L^2+1,\,3\}}$ and
$c_4=\max\{\f{8c_1^2}{c_2^2}+1,\,2\}$.
\end{lemma}
\no {\it Proof}. As $\gamma^{k+1}$ is a minimizer in (\ref{ch4}), we derive using
the estimates (\ref{ch2}), (\ref{es1})-(\ref{es2}), equality \eqref{es7} and the Cauchy-Schwarz inequality
\beqnx
&&\|u'(\gamma^k)(\gamma^{k+1}-\gamma^k)-(z^\delta-u(\gamma^k))\|_{\Gamma_a}^2\leq J(\gamma^{k+1})\leq J(\gamma^*)\nb\\
&=&\|u'(\gamma^k)(\gamma^*-\gamma^k)-(u(\gamma^*)-u(\gamma^k))+u(\gamma^*)-z^\delta\|_{\Gamma_a}^2
+\beta_k\|\gamma^*-\gamma^k\|_{\Gamma_i}^2\nb\\
&\leq&2c_1^2\|\gamma^{k}-\gamma^*\|^4_{\Gamma_i}+2\delta^2+
\|u(\gamma^k)-u(\gamma^*)+u(\gamma^*)-z^\delta\|^2_{\Gamma_a}\|\gamma^*-\gamma^k\|_{\Gamma_i}^2\nb\\
&\leq&2c_1^2\|\gamma^{k}-\gamma^*\|^4_{\Gamma_i}+2\delta^2+2L^2\|\gamma^*-\gamma^k\|_{\Gamma_i}^4+
2\delta^2\|\gamma^*-\gamma^k\|_{\Gamma_i}^2\nb\\
&\leq&(2c_1^2+2L^2+1)\|\gamma^{k}-\gamma^*\|^4_{\Gamma_i}+(2+\delta^2)\delta^2\nb\\
&\leq&\max\{2c_1^2+2L^2+1,\,3\}(\|\gamma^{k}-\gamma^*\|^4_{\Gamma_i}+\delta^2),
\eqnx
which implies \eqref{es11} immediately.

Again, using the minimizing property of $\gamma^{k+1}$ in (\ref{ch4}) and
the estimates (\ref{ch2}) and (\ref{es2}), we can deduce as follows:
\beqn
&&\|\gamma^{k+1}-\gamma^k\|^2_{\Gamma_i}\leq\f{1}{\beta_k}J(\gamma^{k+1})\leq\f{1}{\beta_k}J(\gamma^*)\nb\\
&=&\f{1}{\beta_k}\|u'(\gamma^k)(\gamma^*-\gamma^k)-(z^\delta-u(\gamma^k))\|_{\Gamma_a}^2
+\|\gamma^*-\gamma^k\|_{\Gamma_i}^2\nb\\
&\leq&\f{1}{\beta_k}(2c_1^2\|\gamma^{k}-\gamma^*\|^4_{\Gamma_i}+2\delta^2)+\|\gamma^*-\gamma^k\|_{\Gamma_i}^2\nb\\
&=&\f{2c_1^2\|\gamma^{k}-\gamma^*\|^4_{\Gamma_i}}{\|u(\gamma^k)-z^\delta\|^2_{\Gamma_a}}+
\f{2\delta^2}{\|u(\gamma^k)-z^\delta\|^2_{\Gamma_a}}+\|\gamma^*-\gamma^k\|_{\Gamma_i}^2.
\label{es10}
\eqn
As stated in (\ref{aj1}),
the iterative process (\ref{ch4}) terminates if
$c_2\|\gamma^k-\gamma^*\|_{\Gamma_i}<2\delta$ or $\|u(\gamma^k)-z^\delta\|_{\Gamma_a}<\sqrt{\delta}$. Otherwise
we have
$c_2\|\gamma^k-\gamma^*\|_{\Gamma_i}\geq2\delta$ and $\|u(\gamma^k)-z^\delta\|_{\Gamma_a}\geq\sqrt{\delta}$.
Then we can easily see that $\f{2\delta^2}{\|u(\gamma^k)-z^\delta\|^2_{\Gamma_a}}\leq2\delta$
and
 \beqnx
&&\|u(\gamma^k)-z^\delta\|_{\Gamma_a}\geq
\|u(\gamma^k)-u(\gamma^*)\|_{\Gamma_a}-\|u(\gamma^*)-z^\delta\|_{\Gamma_a}\\
&\geq& c_2\|\gamma^k-\gamma^*\|_{\Gamma_i}-\delta
\geq c_2\|\gamma^k-\gamma^*\|_{\Gamma_i}-\frac{c_2}{2}\|\gamma^k-\gamma^*\|_{\Gamma_i}
=\frac{c_2}{2}\|\gamma^k-\gamma^*\|_{\Gamma_i}\,.
\eqnx
Now the desired result \eqref{es18} follows readily from these two estimates and (\ref{es10}).  $\Box$

\begin{lemma}\label{lem:oestimation1}
Under the conditions (\ref{es8}), (\ref{es7}) and (\ref{aj1}), let $\gamma^k$ and $\gamma^{k+1}$ be
two consequent iterates generated by the iteration (\ref{ch4}) such that
both $\gamma^k$ and $\gamma^{k+1}$ lie in $N(\gamma^*,b)$, then
\beqn
\|\gamma^{k+1}-\gamma^*\|_{\Gamma_i}\leq c_5(\|\gamma^{k}-\gamma^*\|^2_{\Gamma_i}+\delta),
\label{c3}
\eqn
where constant $c_5$ is given explicitly by $c_5={(c_3+c_1c_4+1)}/{c_2}$.
\end{lemma}
\no {\it Proof}. It follows from (\ref{es8}),  (\ref{es2}), (\ref{es11}) and (\ref{es18}) that
\beqnx
&&c_2\|\gamma^{k+1}-\gamma^*\|_{\Gamma_i}\leq\|u(\gamma^{k+1})-u(\gamma^*)\|_{\Gamma_a}\leq
\|u(\gamma^{k+1})-z^\delta\|_{\Gamma_a}+\delta\nb\\
&=&\|u'(\gamma^k)(\gamma^{k+1}-\gamma^k)+u(\gamma^k)-z^\delta-
\{u'(\gamma^k)(\gamma^{k+1}-\gamma^k)+u(\gamma^k)-u(\gamma^{k+1})\}\|_{\Gamma_a}+\delta\nb\\
&\leq&\|u'(\gamma^k)(\gamma^{k+1}-\gamma^k)+u(\gamma^k)-z^\delta)\|_{\Gamma_a}
+\|u'(\gamma^k)(\gamma^{k+1}-\gamma^k)+u(\gamma^k)-u(\gamma^{k+1})\|_{\Gamma_a}+\delta\nb\\
&\leq&c_3(\|\gamma^{k}-\gamma^*\|^2_{\Gamma_i}+\delta)+c_1\|\gamma^{k+1}-\gamma^k\|^2_{\Gamma_i}+\delta\nb\\
&\leq&c_3(\|\gamma^{k}-\gamma^*\|^2_{\Gamma_i}+\delta)+
c_1c_4(\|\gamma^{k}-\gamma^*\|^2_{\Gamma_i}+\delta)+\delta\nb\\
&\leq&(c_3+c_1c_4+1)(\|\gamma^{k}-\gamma^*\|^2_{\Gamma_i}+\delta),
\eqnx
which implies the estimate \eqref{c3}.  $\Box$

In order to establish the quadratic convergence of the L-M iteration, we now emphasize
the dependence of all the constants $c_1$, $\cdots$, $c_5$ in our previous estimates on the radius $b$
of the ball $N(\gamma^*,b)$.
First, we know both constants $c_1$ and $c_3$ in \eqref{es2} and \eqref{es11} are independent of $b$.
But constant $c_2$ in (\ref{es8}) depends on this radius $b$,
so we will write $c_2(b)$ to emphasize this dependence.
Similarly, we can write the constants $c_4$ and $c_5$ in the estimates
\eqref{es18} and \eqref{c3} as $c_4(b)$ and $c_5(b)$.

We are now ready to establish our major convergence results in this work, the quadratic convergence
and quadratic rate of convergence for
the Levenberg-Marquardt iteration (\ref{ch4}). For simplicity, we set
\beqn
&&r(b,\delta)=\min\Big\{b,\,\frac{b-\sqrt{c_4(\frac{1}{3})\delta}}{\sqrt{c_4(\frac{1}{3})}+1}\Big\},
\q \alpha=c_5(\frac{1}{3})\Big(\sqrt{c_4(\frac{1}{3})}+1\Big)\,,\nb\\
&&\beta(\delta)=\delta c_5(\frac{1}{3}) \Big(\sqrt{c_4(\frac{1}{3})}+1\Big) +\sqrt{c_4(\frac{1}{3})\delta}\,. \label{aj2}
\eqn
We can readily see from  (\ref{es1}) and assumption (\ref{es8}) that $c_2(b)\leq L$,
using which we can directly check from the definitions of $c_4$ and $c_5$ that
$c_5(b)>\sqrt{2}$ and $c_4(b)\geq 2$. Using these we know $\alpha >3$, and
${1+\sqrt{1-4\alpha\beta(\delta)}}<(2\alpha)/3$, which will be used in the following
theorem.
\begin{theorem}\label{lem:quad}
Under the conditions (\ref{es8}), (\ref{es7}) and (\ref{aj1}),
we assume $\delta$ is small enough such that $r(b,\delta)>0$, $1-4\alpha\beta(\delta)> 0$
and choose $b\in [\frac{1-\sqrt{1-4\alpha\beta(\delta)}}{2\alpha},\,
\frac{1+\sqrt{1-4\alpha\beta(\delta)}}{2\alpha}]$, then for any $\gamma^0\in N(\gamma^*,r(b,\delta))$,
the sequence $\{\gamma^k\}$ generated by
(\ref{ch4}) stays always in $N(\gamma^*,b)$ and satisfies
\beqnx
\|\gamma^{k+1}-\gamma^*\|_{\Gamma_i}\leq c_5(\frac{1}{3})(\|\gamma^{k}-\gamma^*\|^2_{\Gamma_i}+\delta).
\eqnx
\end{theorem}
\no {\it Proof}.
From the results of Lemma \ref{lem:oestimation1}, we only need to show that the sequence $\{\gamma^k\}$ generated by
(\ref{ch4}) stays always in $N(\gamma^*,b)$.
This is proved below by the mathematical induction.

First, by the choice $b$ we know $0<b<{1}/{3}$, we know
$\gamma^0\in N(\gamma^*,r(b,\delta))\subset N(\gamma^*, {1}/{3})$. Then by the
triangle inequality and estimate (\ref{es18}) with $b={1}/{3}$ we can deduce
\beqnx
\|\gamma^1-\gamma^*\|_{\Gamma_i}&\leq&\|\gamma^1-\gamma^0\|_{\Gamma_i}+\|\gamma^0-\gamma^*\|_{\Gamma_i}
\leq(\sqrt{c_4(\frac{1}{3})}+1)\|\gamma^0-\gamma^*\|_{\Gamma_i}+\sqrt{c_4(\frac{1}{3})\delta}\\
&\leq&(\sqrt{c_4(\frac{1}{3})}+1)r(b,\delta)+\sqrt{c_4(\frac{1}{3})\delta}\leq b,
\eqnx
which implies $\gamma^1\in N(\gamma^*,b)\subset N(\gamma^*, {1}/{3})$.

Now we show $\gamma^{l+1}\in N(\gamma^*,b)$ if
$\gamma^k\in N(\gamma^*,b)$ for $k=1,...,l$.
Indeed, we deduce from
the triangle inequality, the estimate (\ref{es18})
for $\gamma^{l}\in N(\gamma^*,b)$ and the estimate (\ref{c3})
for $\gamma^l,\gamma^{l-1}\in N(\gamma^*,b)$  that

\beqn
\|\gamma^{l+1}-\gamma^*\|_{\Gamma_i}&\leq&\|\gamma^{l+1}-\gamma^l\|_{\Gamma_i}+\|\gamma^l-\gamma^*\|_{\Gamma_i}\nb\\
&\leq&(\sqrt{c_4(\frac{1}{3})}+1)\|\gamma^l-\gamma^*\|_{\Gamma_i}+\sqrt{c_4(\frac{1}{3})\delta}\nb\\
&\leq&(\sqrt{c_4(\frac{1}{3})}+1)c_5(\frac{1}{3})(\|\gamma^{l-1}-\gamma^*\|^2_{\Gamma_i}+\delta)+\sqrt{c_4(\frac{1}{3})\delta}\nb\\
&\leq& \alpha b^2+\beta(\delta)
\leq b\,,
\label{es20}
\eqn
which implies that $\gamma^{l+1}\in N(\gamma^*,b)$, if it holds that
$\alpha b^2+\beta(\delta)\leq b$.

To see $\alpha b^2+\beta(\delta)\leq b$, we define a quadratic functional
$
f(b)=\alpha b^2-b+\beta(\delta)\,.
$
As $1-4\alpha\beta(\delta)> 0$, it is easy to see
$0<b_1=\frac{1-\sqrt{1-4\alpha\beta(\delta)}}{2\alpha}< 1/6$ and
$b_1<b_2=\frac{1+\sqrt{1-4\alpha\beta}}{2\alpha}
<1/3$, and $b_1$ and $b_2$ are two
solutions of $f(b)=0$. Clearly for any $b\in [b_1,\,b_2]$,
we know $f(b)\leq0$, namely, $\alpha b^2+\beta(\delta)\leq b$.
$\Box$

\subsection{Surrogate functional technique}\label{sec:ddmrobin}
In each step of the L-M iteration we have to solve the minimization problem (\ref{ch4}).
Let us now derive its optimality system, i.e., $J'(\gamma^{k+1})\xi=0$ for any $\xi\in L^2(\Gamma_i)$.
By direct computations, we have
\beqnx
J'(\gamma)\xi&=&2\langle u'(\gamma^k)(\gamma-\gamma^k)-(z^\delta-u(\gamma^k)),u'(\gamma^k)(\xi)\rangle_{\Gamma_a}
+2\beta_k\langle\gamma-\gamma^k,\xi\rangle_{\Gamma_i}\\
&=&2\langle u(\gamma^k)\{u'(\gamma^k)^*(\frac{u'(\gamma^k)(\gamma-\gamma^k)-(z^\delta-u(\gamma^k))}{u(\gamma^k)})\},
\xi\rangle_{\Gamma_i}+2\beta_k\langle\gamma-\gamma^k,\xi\rangle_{\Gamma_i},
\eqnx
where we have used the adjoint relation (\ref{ch8}). This is equivalent to the following
equation:
\beqn
u(\gamma^k)\{u'(\gamma^k)^*(\frac{u'(\gamma^k)(\gamma-\gamma^k)}{u(\gamma^k)})\}
+\beta_k(\gamma-\gamma^k)=u(\gamma^k)\{u'(\gamma^k)^*(\frac{z^\delta-u(\gamma^k)}{u(\gamma^k)})\}.
\label{nece1}
\eqn
So we have to solve this rather complicated linear system (whose discretized system is highly ill-conditioned)
to get the solution $\gamma^{k+1}$ at each iteration of (\ref{ch4}), e.g., by some
iterative method.  This is still difficult and computationally very expensive.

Next, we shall make use of the surrogate functional technique
to greatly simplify the solution to the minimization (\ref{ch4}), resulting in an explicit solution
at each iteration. The resultant algorithm is computationally much less expensive.
The surrogate functional technique was studied in \cite{dau04} for solving a linear inverse operator equation of
the form $Kh=f$. We now construct a surrogate functional $J^s(\gamma,\gamma^{k})$ of $J(\gamma)$
in (\ref{ch4}):
\bb
J^s(\gamma,\gamma^{k})=J(\gamma)
+A\|\gamma-\gamma^{k}\|_{\Gamma_i}^2-\|u'(\gamma^{k})(\gamma-\gamma^{k})\|_{\Gamma_a}^2,
\label{che1}
\ee
where $A$ can be any positive constant such that $\|u'(\gamma^{k})d\|_{\Gamma_a}^2\le A\|d\|_{\Gamma_i}^2$
for all $d\in L^2(\Gamma_i)$. Next, we will simplify the expression (\ref{che1}).
Using the adjoint relation (\ref{ch8}), we can rewrite $J^s(\gamma,\gamma^{k})$ as follows:
\beqn
J^s(\gamma,\gamma^{k})&=&
\|u'(\gamma^k)(\gamma-\gamma^k)-(z^\delta-u(\gamma^k))\|_{\Gamma_a}^2
+\beta_k\|\gamma-\gamma^k\|_{\Gamma_i}^2+A\|\gamma-\gamma^{k}\|_{\Gamma_i}^2\nb\\
&&-\|u'(\gamma^{k})(\gamma-\gamma^{k})\|_{\Gamma_a}^2\nb\\
&=&-2\langle u'(\gamma^k)(\gamma-\gamma^k),\,
z^\delta-u(\gamma^k)\rangle_{\Gamma_a}+\beta_k\|\gamma-\gamma^k\|_{\Gamma_i}^2
+A\|\gamma-\gamma^{k}\|_{\Gamma_i}^2\nb\\
&&+\|z^\delta-u(\gamma^k)\|^2_{\Gamma_a}\nb\\
&=&-2\langle \gamma-\gamma^k,\,u(\gamma^k)
\{u'(\gamma^k)^*(\frac{z^\delta-u(\gamma^k}{u(\gamma^k)})\}\rangle_{\Gamma_i}
+\beta_k\|\gamma-\gamma^k\|_{\Gamma_i}^2\nb\\
&&+A\|\gamma-\gamma^{k}\|_{\Gamma_i}^2
+\|z^\delta-u(\gamma^k)\|^2_{\Gamma_a}\nb\\
&=&A\|\gamma-\gamma^{k}-\f{1}{A}u(\gamma^k)
\{u'(\gamma^k)^*(\frac{z^\delta-u(\gamma^k)}{u(\gamma^k)})\}\|_{\Gamma_i}^2
+\beta_k\|\gamma-\gamma^k\|_{\Gamma_i}^2\nb\\
&&+\Big\{ \|z^\delta-u(\gamma^k)\|^2_{\Gamma_a}-A\|\f{1}{A}u(\gamma^k)
\{u'(\gamma^k)^*(\frac{z^\delta-u(\gamma^k)}{u(\gamma^k)})\}\|_{\Gamma_i}^2\Big\}.
\label{che2}
\eqn
We can see that the last term above is independent of $\gamma$, so does not affect the minimization.
Hence we will drop that term in the functional $J^s(\gamma,\gamma^{k})$ and obtain
\beqn
\min_{\gamma\in K}J^s(\gamma,\gamma^{k})
=\min_{\gamma\in K}A\Big\|\gamma-\gamma^{k}-\f{1}{A}u(\gamma^k)
\Big\{u'(\gamma^k)^*(\frac{z^\delta-u(\gamma^k)}{u(\gamma^k)})\Big\}\Big\|_{\Gamma_i}^2
+\beta_k\|\gamma-\gamma^k\|_{\Gamma_i}^2.
\label{che3}
\eqn
This is a simple quadratic minimization, and we can compute its minimizer exactly:
\beqn
\underset{\gamma\in K}{\rm argmin}\,J^s(\gamma,\gamma^{k})
=\gamma^k+\frac{1}{A+\beta_k}u(\gamma^k)
\Big\{u'(\gamma^k)^*(\frac{z^\delta-u(\gamma^k)}{u(\gamma^k)})\Big\}.
\label{che5}
\eqn
This motivates us with the following reconstruction algorithm for the Robin coefficient
in (\ref{ch1}), which is clearly much easier and computationally much less expensive than
solving the minimization (\ref{ch4}) directly.

\begin{algorithm}\label{al:ell}

 Choose a tolerance parameter $\epsilon>0$ and an initial value
$\gamma^{0}$, and set $k:=0$.
\bn
\item Compute $\gamma^{k+1}$:
\beqn
\gamma^{k+1}&=&\gamma^k+\frac{1}{A+\beta_k}u(\gamma^k)
\Big\{u'(\gamma^k)^*(\frac{z^\delta-u(\gamma^k)}{u(\gamma^k)})\Big\}.
\label{cc6}
\eqn

\item If $\frac{\|\gamma^{k+1}-\gamma^{k}\|_{\Gamma_i}}{\|\gamma^{k}\|_{\Gamma_i}}
\leq\epsilon$, stop the iteration;
otherwise set $k:=k+1$, go to Step 1.

\en
\end{algorithm}


\section{Parabolic inverse Robin problem and its L-M solution}\label{sec:formulationt}
\setcounter{equation}{0}
\setcounter{figure}{0}
\setcounter{table}{0}

\subsection{Tikhonov regularization for the parabolic inverse Robin problem}
\label{sub:b2}
In this section we first formulate the Levenberg-Marquardt method for solving
the nonlinear non-convex optimisation problems resulting from the least-squares formulation
of the parabolic inverse Robin problem as stated in Section\,\ref{sec:intro}, incorporated with Tikhonov regularization
to handle its ill-posedness and instability due to the presence of the noise in the observation data \cite{jin12}.
We assume the noise level in the observation data $z^\delta$ of the true solution $u$
to the parabolic system (\ref{ch1t})
is of order $\delta$, namely
\bb
\int_{0}^T\|u(\gamma^*)-z^\delta\|^2_{\Gamma_a}dt\leq \delta^2.
\label{cch2t}
\ee
where $\gamma^*$ is the true Robin coefficient in the system (\ref{ch1t}).
The parabolic inverse Robin problem is frequently transformed
into the following stabilized minimization system with Tikhonov regularization:
\beqn
\min_{\gamma\in K}{\cal J}(\gamma)
&=& \int_{0}^T\|u(\gamma)-z^\delta\|_{\Gamma_a}^2dt+\beta\|\gamma\|_{\Gamma_i}^2\,.
\label{ch3t}
\eqn
where $\beta$ is the regularization parameter.
The formulation (\ref{ch3t}) was shown to be stable in the sense that its minimizer depends continuously on the change of
the noise in the data $z^\delta$ \cite{jin12}.

%
For the subsequent analysis on the convergence of the Levenberg-Marquardt method
for solving the optimisation (\ref{ch3t}), we shall frequently need the Fr\'{e}chet derivative of the forward
solution $u(\gamma)$ of system (\ref{ch1t}). Let $w:=u'(\gamma)d$ be the Fr\'{e}chet derivative
at direction $d$, then $w\in L^2(0,T;H^1(\Om))$ solves the following system:
\bb \left\{ \begin{array} {rclll}
\p_t w-\nabla\cdot(a(\x)\nabla w) &=&0 & \m{in} &\Om\times (0,T)\,, \\
a(\x)\frac{\p w}{\p n}+\gamma \,w&=&-d\,u(\gamma) & \m{on} &\Gamma_i\times (0,T)\,, \\
a(\x)\frac{\p w}{\p n} &=&0 & \m{on} &\Gamma_a\times (0,T)\,,\\
w(\x,0) &=&0 & \m{in} &\Om\,,
\end{array}
\right. \label{ch6t}
\ee
Let $u'(\gamma)^*$ be the adjoint operator of the Fr\'{e}chet derivative $u'(\gamma)$, then it it easy to verify that
$w^*:=u'(\gamma)^*p\in  L^2(0,T;H^1(\Om))$ at a general direction $p$
solves the following system:
\bb \left\{ \begin{array} {rclll}
-\p_t w^*-\nabla\cdot(a(\x)\nabla w^*) &=&0 & \m{in} &\Om\times [0,T]\,, \\
a(\x)\frac{\p w^*}{\p n}+\gamma \,w^*&=&0 & \m{on} &\Gamma_i\times [0,T]\,, \\
a(\x)\frac{\p w^*}{\p n} &=&-p\,u(\gamma) & \m{on} &\Gamma_a\times [0,T]\,,\\
w^*(\x,T) &=&0 & \m{in} &\Om\,.
\end{array}
\right. \label{ch7t}
\ee

The following lemma gives an important relation for our later analysis.
\begin{lemma}\label{lem:adjointt}
It holds for any directions $d$ and $p$ that
\bb
\int_{0}^T\langle u(\gamma)w,\,p\rangle_{\Gamma_a}dt=\langle
d,\,\int_0^Tu(\gamma)w^*dt\rangle_{\Gamma_i}.
\label{ch8t}
\ee
\end{lemma}
\no {\it Proof}. For any $\varphi, \psi\in L^2(0,T;H^1(\Om))$, we can readily derive the
variational forms of (\ref{ch6t})
and (\ref{ch7t}):
\beqn
\int_0^T\int_\Om \p_tw\varphi d\x dt&+&
\int_0^T\int_\Om a(\x)\nabla w\cdot\nabla \varphi d\x dt
=
\int_0^T\int_{\Gamma_i}(-du(\gamma)-\gamma w)\varphi ds dt,
\label{ch9t}
\eqn
\beqn
-\int_0^T\int_\Om \p_tw^*\psi d\x dt&+&
\int_0^T\int_\Om a(\x)\nabla
w^*\cdot\nabla \psi d\x dt
\nb\\&=&
-\int_0^T\int_{\Gamma_i}\gamma w^*\psi ds dt
-\int_{0}^T\int_{\Gamma_a}p u(\gamma)\psi ds dt.
\label{ch10t}
\eqn
By integrating by parts with respect to $t$ in (\ref{ch10t}), we derive
\beqn
&&-\int_\Om w^*(\x,T)\psi(\x,T) d\x+\int_\Om w^*(\x,0)\psi(\x,0) d\x
+\int_0^T\int_\Om \p_t\psi w^* d\x dt\nb\\
&&+\int_0^T\int_\Om a(\x)\nabla
w^*\cdot\nabla \psi d\x dt=
-\int_0^T\int_{\Gamma_i}\gamma w^*\psi ds dt
-\int_{0}^T\int_{\Gamma_a}p u(\gamma)\psi ds dt.
\label{ch11t}
\eqn
Taking $\varphi=w^*$ in (\ref{ch9t}), $\psi=w$ in (\ref{ch11t})
and noting that $w^*(\x,T)=w(\x,0)=0$, we can easily find that
\beqnx
\int_{0}^T\langle u(\gamma)w,\,p\rangle_{\Gamma_a}dt=\int_0^T\langle
d,\,u(\gamma)w^*\rangle_{\Gamma_i}dt=\langle
d,\,\int_0^Tu(\gamma)w^* dt\rangle_{\Gamma_i}.~~~~~~~~\Box
\eqnx

\subsection{Levenberg-Marquardt method and its convergence}\label{sec:linearp}
The nonlinearity of the parabolic forward solution $u(\gamma)$ of the system (\ref{ch1t})
makes the minimization (\ref{ch3t}) highly nonlinear and non-convex
with respect to the Robin coefficient $\gamma$, as well as strongly unstable
at discrete level with fine mesh sizes and time step sizes due to the severe ill-posedness
of the inverse problem and the fact that noise is always present in the observation data.
To alleviate these difficulties in numerical solutions, we propose to solve the minimization (\ref{ch3t})
by the Levenberg-Marquardt iteration:
\beqn
J(\gamma^{k+1})=\min_{\gamma\in K}J(\gamma)
=: \int_{0}^T\|u'(\gamma^k)(\gamma-\gamma^k)-(z^\delta-u(\gamma^k))\|_{\Gamma_a}^2dt
+\beta_k\|\gamma-\gamma^k\|_{\Gamma_i}^2\,.
\label{ch4t}
\eqn
Before our study of the convergence of the iteration (\ref{ch4t}),
we establish some important auxiliary results.
\begin{lemma}\label{lem:estimationp}
Let  $u(\gamma)$ be the forward operator of the system (\ref{ch1t}) and
$u(\gamma)\in L^2(0,T;L^{\infty}(\Gamma_i))$  for $\gamma\in K$, then
there exist positive constants $\bar L$ and $\bar c_1$ such that the following estimates hold
for any $\gamma,\,\bar\gamma\in K$:
%
\beqn
\int_{0}^T\|u(\gamma)-u(\bar \gamma)\|^2_{\Gamma_a}dt&\leq&\bar L\|\gamma-\bar \gamma\|^2_{\Gamma_i},\label{es1t}\\
\int_{0}^T\|u'(\bar \gamma)(\gamma-\bar \gamma)-(u(\gamma)-u(\bar \gamma))\|^2_{\Gamma_a}dt
&\leq&\bar c_1\|\gamma-\bar \gamma\|^4_{\Gamma_i}\,.\label{es2t}
\eqn
\end{lemma}
\no {\it Proof}. From the variational form of the system (\ref{ch1t}), we can easily see
for any $\varphi\in L^2(0,T;H^1(\Om))$ that
\beqn
&&\int_\Om \p_t(u(\gamma)-u(\bar \gamma))\varphi d\x+
\int_\Om a(\x)\nabla (u(\gamma)-u(\bar \gamma))\cdot\nabla \varphi d\x d\x\nb\\
&+&\int_{\Gamma_i}\bar\gamma (u(\gamma)-u(\bar \gamma)) \varphi ds=
-\int_{\Gamma_i}(\gamma-\bar\gamma)u(\gamma)\varphi ds.
\label{es3t}
\eqn
Taking $\varphi=u(\gamma)-u(\bar \gamma)$ in (\ref{es3t}) and integrating  by parts with respect
to $t$ over $[0,\tau]$ for $\tau\in [0,T]$, then using the Cauchy-Schwarz inequality, we derive
\beqnx
&&\frac{1}{2}\|u(\gamma)(\x,\tau)-u(\bar \gamma)(\x,\tau)\|^2_{\Om}
+\underline a\int_0^\tau\|\nabla u(\gamma)-\nabla u(\bar \gamma)\|^2_{\Om}dt\\
&&+\gamma_1\int_0^\tau\|u(\gamma)-u(\bar \gamma)\|^2_{\Gamma_i}dt
\leq|\int_0^\tau\int_{\Gamma_i}(\gamma-\bar\gamma)u(\gamma)(u(\gamma)-u(\bar \gamma)) dsdt|\\
&&\leq\|u(\gamma)\|_{L^2(0,T;L^{\infty}(\Gamma_i))}\sqrt{\tau}\|\gamma-\bar\gamma\|_{\Gamma_i}
\|u(\gamma)-u(\bar \gamma)\|_{L^2(0,\tau;L^{2}(\Gamma_i))}\\
&&\leq C\|\gamma-\bar\gamma\|_{\Gamma_i}
\|u(\gamma)-u(\bar \gamma)\|_{L^2(0,\tau;L^{2}(\Gamma_i))}.
\eqnx
Now a direct application of the Young's inequality gives
\beqnx
\|u(\gamma)(\x,\tau)-u(\bar \gamma)(\x,\tau)\|^2_{\Om}
+\int_0^\tau\|\nabla u(\gamma)-\nabla u(\bar \gamma)\|^2_{\Om}dt
\leq C\|\gamma-\bar\gamma\|^2_{\Gamma_i},
\eqnx
from which and the trace theorem, we obtain
\beqnx
\int_0^T\|u(\gamma)-u(\bar \gamma)\|^2_{\Gamma_a}dt
\leq C\int_0^T\|u(\gamma)-u(\bar \gamma)\|^2_{1,\Om}dt
\leq \bar L \|\gamma-\bar\gamma\|^2_{\Gamma_i}.
\eqnx
To verify the estimate (\ref{es2t}), we first have by taking $\varphi=u'(\gamma)d$ in (\ref{ch9t})
and then following the same technique as we did for (\ref{es1t}) that
\bb
\int_0^T\|u'(\gamma)d\|^2_{\p\Om}dt\leq C\|\gamma-\bar\gamma\|^2_{\Gamma_i}\,.
\label{es5t}
\ee
Next, we take $\gamma=\bar\gamma$ and $d=\gamma-\bar\gamma$ in (\ref{ch9t}) to deduce
\beqn
\int_0^T\int_\Om \p_t(u'(\bar\gamma)(\gamma-\bar\gamma))\varphi d\x dt&+&
\int_0^T\int_\Om a(\x)\nabla (u'(\bar\gamma)(\gamma-\bar\gamma))\cdot\nabla \varphi d\x dt\nb\\
&=&
\int_0^T\int_{\Gamma_i}(-(\gamma-\bar\gamma)u(\bar\gamma)-
\bar\gamma (u'(\bar\gamma)(\gamma-\bar\gamma)))\varphi ds dt.
\label{es4t}
\eqn
Subtracting (\ref{es4t}) from (\ref{es3t}) gives
\beqnx
&&\int_0^T\int_\Om \p_t(u'(\bar\gamma)(\gamma-\bar\gamma)-(u(\gamma)-u(\bar \gamma)))\varphi d\x dt+
\int_0^T\int_\Om a(\x)\nabla (u'(\bar\gamma)(\gamma-\bar\gamma)-(u(\gamma)-u(\bar \gamma)))\cdot\nabla \varphi d\x dt\nb\\
&&+\int_0^T\int_{\Gamma_i}
\bar\gamma (u'(\bar\gamma)(\gamma-\bar\gamma)-(u(\gamma)-u(\bar \gamma)))\varphi ds dt
=\int_0^T\int_{\Gamma_i}(\gamma-\bar\gamma)(u(\gamma)-u(\bar \gamma))\varphi ds dt.
\eqnx
Now applying the trace theorem, Lagrange mean value theorem and estimate (\ref{es5t}), we can derive
\beqnx
&&\int_{0}^T\|u'(\bar \gamma)(\gamma-\bar \gamma)-(u(\gamma)-u(\bar \gamma))\|^2_{\Gamma_a}dt
\leq C\int_0^T\|u'(\bar \gamma)(\gamma-\bar \gamma)-(u(\gamma)-u(\bar \gamma))\|^2_{1,\Om}dt\\
&\leq&C\int_0^T \|(\gamma-\bar\gamma)(u(\gamma)-u(\bar \gamma))\|^2_{\Gamma_i}dt
=C\int_0^T \|(\gamma-\bar\gamma) u'(\xi)(\gamma-\bar\gamma)\|^2_{\Gamma_i}dt
\leq \bar c_1\|\gamma-\bar\gamma\|^4_{\Gamma_i},
\eqnx
where $\xi$ is some element in $K$ between $\gamma$ and $\bar\gamma$.   $\Box$

Now we are ready to establish a quadratic rate on the convergence of the Levenberg-Marquardt method,
under the following basic condition:
\bb
\int_{0}^T\|u(\gamma)-u(\gamma^*)\|^2_{\Gamma_a}dt\geq \bar c_2\|\gamma-\gamma^*\|^2_{\Gamma_i},
 \q \forall\,\gamma\in N(\gamma^*,b)
\label{es8t}
\ee
where $\bar c_2$ and $b$ are two positive constants with $b\in (0,1)$.
Assumption (\ref{es8t}) is
the frequently adopted basic condition to ensure the quadratic convergence of the Levenberg-Marquardt
method for most direct nonlinear optimization problems [4] [15], so it is natural to bring it to the current
nonlinear ill-posed parabolic inverse problem. The condition (\ref{es8t}) may be viewed as a direct motivation
of the local Lipschitz stability of the parabolic inverse Robin problems (see estimate (\ref{ccc6}) in Theorem \ref{thm:localp}).

It is a well-known technical difficulty in a practical numerical realisation of any Tikhonov regularised optimisation system like
the ones (\ref{ch3t}) and \eqref{ch4t} to choose a reasonable and effective regularization parameter
$\beta$ or $\beta_k$.
The same as we did in Section\,\ref{sec:linear}, one of the important novelties of this work is our suggestion of
a very simple and easy implementable choice of the parameter $\beta_k$ based on the rule:
\bb
\beta_k=\int_{0}^T\|u(\gamma^k)-z^\delta\|^2_{\Gamma_a}dt,
\label{es7t}
\ee
And surprisingly, as we shall demonstrate below, this choice of the regularization parameter $\beta_k$ still ensures
a quadratical convergence of the resulting Levenberg-Marquardt iteration (\ref{ch4t}).

Considering the presence of the noise (see \eqref{cch2t}), it is reasonable for us to terminate the L-M iteration (\ref{ch4t})
when its minimizer $\gamma^k$ is accurate enough in terms of the noise level, more specifically,
we shall terminate the iteration if the following criterion is realised:
\bb
\sqrt{\bar c_2}\|\gamma^k-\gamma^*\|_{\Gamma_i}<2\delta \q \m{or} \q
\int_{0}^T\|u(\gamma^k)-z^\delta\|^2_{\Gamma_a}dt<\delta\,.
\label{aj1t}
\ee

\begin{lemma}\label{lem:oestimationp}
Under the conditions (\ref{es8t}), (\ref{es7t}) and (\ref{aj1t}),
if $\gamma^k\in N(\gamma^*,b)$, then $\gamma^{k+1}$ generated by the iteration (\ref{ch4t}) satisfies
\beqn
\int_{0}^T\|u'(\gamma^k)(\gamma^{k+1}-\gamma^k)-(z^\delta-u(\gamma^k))\|^2_{\Gamma_a}dt&\leq&
\bar c_3(\|\gamma^{k}-\gamma^*\|^4_{\Gamma_i}+\delta^2), \label{es11t}\\
\|\gamma^{k+1}-\gamma^k\|^2_{\Gamma_i}
&\leq& \bar c_4(\|\gamma^{k}-\gamma^*\|^2_{\Gamma_i}+\delta),
\label{es18t}
\eqn
where $\bar c_3$ and $\bar c_4$ are two constants given explicitly by
$\bar c_3=\max\{2\bar c_1+2\bar L+1,\,3\}$ and
$\bar c_4=\max\{\f{8\bar c_1}{\bar c_2}+1,\,2\}$.
\end{lemma}
\no {\it Proof}. As $\gamma^{k+1}$ is a minimizer in (\ref{ch4t}), we can derive by using
the estimates (\ref{cch2t}), (\ref{es1t})-(\ref{es2t}), the equality (\ref{es7t}) and Cauchy-Schwarz inequality that
\beqnx
&&\int_{0}^T\|u'(\gamma^k)(\gamma^{k+1}-\gamma^k)-(z^\delta-u(\gamma^k))\|_{\Gamma_a}^2dt
\leq J(\gamma^{k+1})\leq J(\gamma^*)\nb\\
&=&\int_{0}^T\|u'(\gamma^k)(\gamma^*-\gamma^k)-(u(\gamma^*)-u(\gamma^k))+u(\gamma^*)-z^\delta\|_{\Gamma_a}^2dt
+\beta_k\|\gamma^*-\gamma^k\|_{\Gamma_i}^2\nb\\
&\leq&2\bar c_1\|\gamma^{k}-\gamma^*\|^4_{\Gamma_i}+2\delta^2+
\int_{0}^T\|u(\gamma^k)-u(\gamma^*)+u(\gamma^*)-z^\delta\|^2_{\Gamma_a}dt\|\gamma^*-\gamma^k\|_{\Gamma_i}^2\nb\\
&\leq&2\bar c_1\|\gamma^{k}-\gamma^*\|^4_{\Gamma_i}+2\delta^2+2\bar L\|\gamma^*-\gamma^k\|_{\Gamma_i}^4+
2\delta^2\|\gamma^*-\gamma^k\|_{\Gamma_i}^2\nb\\
&\leq&(2\bar c_1+2\bar L+1)\|\gamma^{k}-\gamma^*\|^4_{\Gamma_i}+(2+\delta^2)\delta^2\\
&\leq&
\max\{2\bar c_1+2\bar L+1,\,3\}(\|\gamma^{k}-\gamma^*\|^4_{\Gamma_i}+\delta^2)\\
&\equiv&
\bar c_3(\|\gamma^{k}-\gamma^*\|^4_{\Gamma_i}+\delta^2).
\eqnx

Again, using the minimizing property of $\gamma^{k+1}$ in (\ref{ch4t}) and
the estimates (\ref{cch2t}) and (\ref{es2t}), we can deduce as follows:
\beqn
&&\|\gamma^{k+1}-\gamma^k\|^2_{\Gamma_i}\leq\f{1}{\beta_k}J(\gamma^{k+1})\leq\f{1}{\beta_k}J(\gamma^*)\nb\\
&=&\f{1}{\beta_k}\int_{0}^T\|u'(\gamma^k)(\gamma^*-\gamma^k)-(z^\delta-u(\gamma^k))\|_{\Gamma_a}^2dt
+\|\gamma^*-\gamma^k\|_{\Gamma_i}^2\nb\\
&\leq&\f{1}{\beta_k}(2\bar c_1\|\gamma^{k}-\gamma^*\|^4_{\Gamma_i}+2\delta^2)+\|\gamma^*-\gamma^k\|_{\Gamma_i}^2\nb\\
&=&\f{2\bar c_1\|\gamma^{k}-\gamma^*\|^4_{\Gamma_i}}{\int_{0}^T\|u(\gamma^k)-z^\delta\|^2_{\Gamma_a}dt}+
\f{2\delta^2}{\int_{0}^T\|u(\gamma^k)-z^\delta\|^2_{\Gamma_a}dt}+\|\gamma^*-\gamma^k\|_{\Gamma_i}^2.
\label{es10t}
\eqn
As stated in (\ref{aj1t}),
the iterative process (\ref{ch4t}) terminates if
$\sqrt{\bar c_2}\|\gamma^k-\gamma^*\|_{\Gamma_i}<2\delta$ or $\int_{0}^T\|u(\gamma^k)-z^\delta\|_{\Gamma_a}^2dt<\delta$.
Otherwise we have
$\sqrt{\bar c_2}\|\gamma^k-\gamma^*\|_{\Gamma_i}\geq2\delta$ and $\int_{0}^T\|u(\gamma^k)-z^\delta\|_{\Gamma_a}^2dt\geq\delta$,
which implies
\beqnx
&&(\int_{0}^T\|u(\gamma^k)-z^\delta\|^2_{\Gamma_a}dt)^{\frac{1}{2}}\geq
(\int_{0}^T\|u(\gamma^k)-u(\gamma^*)\|_{\Gamma_a}^2dt)^{\frac{1}{2}}
-(\int_{0}^T\|u(\gamma^*)-z^\delta\|_{\Gamma_a}^2dt)^{\frac{1}{2}}\\
&\geq& \sqrt{\bar c_2}\|\gamma^k-\gamma^*\|_{\Gamma_i}-\delta
\geq \sqrt{\bar c_2}\|\gamma^k-\gamma^*\|_{\Gamma_i}-\frac{\sqrt{\bar c_2}}{2}\|\gamma^k-\gamma^*\|_{\Gamma_i}
=\frac{\sqrt{\bar c_2}}{2}\|\gamma^k-\gamma^*\|_{\Gamma_i}.
\eqnx
Now the desired estimate \eqref{es18t} follows directly from (\ref{es10t}) and the above two estimates.
$\Box$

\begin{lemma}\label{lem:oestimation1p}
Under the conditions (\ref{es8t}), (\ref{es7t}) and (\ref{aj1t}), let $\gamma^k$ and $\gamma^{k+1}$ be
two consequent iterates generated by the iteration (\ref{ch4t}) and satisfy
that $\gamma^k, \, \gamma^{k+1}\in N(\gamma^*,b)$, then
\beqn
\|\gamma^{k+1}-\gamma^*\|_{\Gamma_i}\leq \bar c_5(\|\gamma^{k}-\gamma^*\|^2_{\Gamma_i}+\delta),
\label{c3t}
\eqn
where $\bar c_5$ is a constant given explicitly by $\bar c_5=\sqrt{(2\bar c_3+4\bar c_1\bar c_4^2+1)/{\bar c_2}}$.
\end{lemma}
\no {\it Proof}. By direct computing and triangle inequality it follows from (\ref{es8t}),  (\ref{es2t}), (\ref{es11t}) and (\ref{es18t}) that
\beqnx
&&\bar c_2\|\gamma^{k+1}-\gamma^*\|_{\Gamma_i}^2\leq\int_{0}^T\|u(\gamma^{k+1})-u(\gamma^*)\|^2_{\Gamma_a}dt\leq
\int_{0}^T\|u(\gamma^{k+1})-z^\delta\|^2_{\Gamma_a}dt+\delta^2\nb\\
&=&\int_{0}^T\|u'(\gamma^k)(\gamma^{k+1}-\gamma^k)+u(\gamma^k)-z^\delta-
\{u'(\gamma^k)(\gamma^{k+1}-\gamma^k)+u(\gamma^k)-u(\gamma^{k+1})\}\|_{\Gamma_a}^2dt+\delta^2\nb\\
&\leq&2\int_{0}^T\|u'(\gamma^k)(\gamma^{k+1}-\gamma^k)+u(\gamma^k)-z^\delta)\|_{\Gamma_a}^2dt
+2\bar c_1\|\gamma^{k+1}-\gamma^k\|^4_{\Gamma_i}+\delta^2\nb\\
&\leq&2\bar c_3(\|\gamma^{k}-\gamma^*\|^4_{\Gamma_i}+\delta^2)+
2\bar c_1\bar c_4^2(\|\gamma^{k}-\gamma^*\|^2_{\Gamma_i}+\delta)^2+\delta^2\nb\\
&\leq&(2\bar c_3+4\bar c_1\bar c_4^2+1)(\|\gamma^{k}-\gamma^*\|^4_{\Gamma_i}+\delta^2),
\eqnx
which leads to \eqref{c3t} directly. $\Box$

In order to establish the quadratic convergence of the L-M iteration, we now emphasise
the dependence of all the constants $\bar c_1$, $\cdots$, $\bar c_5$ in our previous estimates on the radius $b$
of the ball $N(\gamma^*,b)$. First, we can easily see that both constants $\bar c_1$ and $\bar c_3$
in \eqref{es2t} and \eqref{es11t} are independent of $b$. But the constant $\bar c_2$ in (\ref{es8t})
depends on this radius $b$, so we will write $\bar c_2(b)$ to emphasize this dependence.
Similarly, we may also write the constants $\bar c_4$ and $\bar c_5$ in the estimates
\eqref{es18t} and \eqref{c3t} as $\bar c_4(b)$ and $\bar c_5(b)$.

We are now ready to establish our major convergence results in this work, the quadratic convergence
and quadratic rate of convergence
for the Levenberg-Marquardt iteration (\ref{ch4t}). For simplicity, we set
%
%
\beqn
&&\bar r(b,\delta)=\min\Big\{b,\,\frac{b-\sqrt{\bar c_4(\frac{1}{3})\delta}}{\sqrt{\bar c_4(\frac{1}{3})}+1}\Big\},
\q \bar \alpha=(\sqrt{\bar c_4(\frac{1}{3})}+1)\bar c_5(\frac{1}{3})\,,\nb\\
&&\bar \beta(\delta)=(\sqrt{\bar c_4(\frac{1}{3})}+1)\bar c_5(\frac{1}{3})\delta+\sqrt{\bar c_4(\frac{1}{3})\delta}\,. \label{aj2}
\eqn
We can directly verify from the definitions of the constants $\bar c_4$ and $\bar c_5$ that
$\bar c_4(b)\geq 2$, $\bar c_5(b)> 2$ and $\bar \alpha >3$.

Following the same arguments as the ones for Theorem\,~\ref{lem:quad}, we can derive

\begin{theorem}\label{lem:quadp}
Under the conditions (\ref{es8t}), (\ref{es7t}) and (\ref{aj1t}),
we assume $\delta$ is small enough such that $\bar r(b,\delta)>0$, $1-4\bar \alpha\bar \beta(\delta)> 0$
and choose $b\in [\frac{1-\sqrt{1-4\bar \alpha\bar \beta(\delta)}}{2\bar \alpha},\,
\frac{1+\sqrt{1-4\bar \alpha\bar \beta(\delta)}}{2\bar \alpha}]$, then for any $\gamma^0\in N(\gamma^*,\bar r(b,\delta))$,
the sequence $\{\gamma^k\}$ generated by
(\ref{ch4}) stay always in $N(\gamma^*,b)$, and
\beqnx
\|\gamma^{k+1}-\gamma^*\|_{\Gamma_i}\leq \bar c_5(\frac{1}{3})(\|\gamma^{k}-\gamma^*\|^2_{\Gamma_i}+\delta).
\eqnx
\end{theorem}

\subsection{Surrogate functional method}\label{sec:ddmrobint}
In each step of the L-M iteration we have to solve the minimization problem (\ref{ch4t}).
Let us now derive its optimality system, i.e., $J'(\gamma^{k+1})\xi=0$ for any $\xi\in L^2(\Gamma_i)$.
By direct computations, we have
%
\beqnx
J'(\gamma)\xi&=&2\int_{0}^T\langle u'(\gamma^k)(\gamma-\gamma^k)-
(z^\delta-u(\gamma^k)),u'(\gamma^k)(\xi)\rangle_{\Gamma_a}dt
+2\beta_k\langle\gamma-\gamma^k,\xi\rangle_{\Gamma_i}\\
&=&2\langle \int_{0}^Tu(\gamma^k)\{u'(\gamma^k)^*(\frac{u'(\gamma^k)(\gamma-\gamma^k)-(z^\delta-u(\gamma^k))}{u(\gamma^k)})\}dt,
\xi\rangle_{\Gamma_i}+2\beta_k\langle\gamma-\gamma^k,\xi\rangle_{\Gamma_i},
\eqnx
where we have used the adjoint relation (\ref{ch8t}). This is equivalent to the following equation:
\beqn
\int_{0}^Tu(\gamma^k)\{u'(\gamma^k)^*(\frac{u'(\gamma^k)(\gamma-\gamma^k)}{u(\gamma^k)})\}dt
+\beta_k(\gamma-\gamma^k)=\int_{0}^Tu(\gamma^k)\{u'(\gamma^k)^*(\frac{z^\delta-u(\gamma^k)}{u(\gamma^k)})\}dt.
\label{nece2}
\eqn
So we have to solve this rather complicated linear system (whose discretized system is strongly ill-conditioned)
to get the solution $\gamma^{k+1}$ at each step of the iteration (\ref{ch4t}), e.g., by some
iterative method. This is still difficult and computationally rather expensive.

Next, we shall make use of the surrogate functional technique
in an attempt to greatly simplify the solution to the minimization (\ref{ch4t}), resulting in an explicit solution
at each iteration. The resultant algorithm is computationally much less expensive.
For the purpose,
we construct an auxiliary surrogate functional $J^s(\gamma,\gamma^{k})$ of $J(\gamma)$ of the form:
\bb
J^s(\gamma,\gamma^{k})=J(\gamma)
+A\|\gamma-\gamma^{k}\|_{\Gamma_i}^2-\int_{0}^T\|u'(\gamma^{k})(\gamma-\gamma^{k})\|_{\Gamma_a}^2dt,
\label{che1p}
\ee
where $A$ is any positive constant such that $\int_{T-\sigma}^T\|u'(\gamma^{k})d\|_{\Gamma_a}^2dt \le
A\|d\|_{\Gamma_i}^2$ for $d\in L^2(\Gamma_i)$. Now we shall convert the functional $J^s(\gamma,\gamma^{k})$
in (\ref{che1p}) into a more explicit representation.
Using the adjoint relation (\ref{ch8t}),
we can rewrite $J^s(\gamma,\gamma^{k})$ as follows:
\beqn
 J^s(\gamma,\gamma^{k})&=&
\int_{0}^T\|u'(\gamma^{k})(\gamma-\gamma^{k})-(z^\delta-u(\gamma^{k}))\|_{\Gamma_a}^2dt
+\beta_k\|\gamma-\gamma^{k}\|_{\Gamma_i}^2+A\|\gamma-\gamma^{k}\|_{\Gamma_i}^2\nb\\
&&-\int_{0}^T\|u'(\gamma^{k})(\gamma-\gamma^{k})\|_{\Gamma_a}^2dt\nb\\
&=&-2\int_{0}^T\langle u'(\gamma^k)(\gamma-\gamma^k),\,
z^\delta-u(\gamma^k)\rangle_{\Gamma_a}dt+\beta_k\|\gamma-\gamma^k\|_{\Gamma_i}^2
+A\|\gamma-\gamma^{k}\|_{\Gamma_i}^2\nb\\
&&+\int_{0}^T\|z^\delta-u(\gamma^{k})\|_{\Gamma_a}^2dt\nb\\
&=&-2\langle \gamma-\gamma^k,\,\int_0^Tu(\gamma^k)
\{u'(\gamma^k)^*(\frac{z^\delta-u(\gamma^k)}{u(\gamma^k)})\}dt\rangle_{\Gamma_i}+\beta_k
\|\gamma-\gamma^k\|_{\Gamma_i}^2\nb\\
&&+A\|\gamma-\gamma^{k}\|_{\Gamma_i}^2
+\int_{0}^T\|z^\delta-u(\gamma^{k})\|_{\Gamma_a}^2dt\nb\\
&=&A\|\gamma-\gamma^{k}-\frac{1}{A}\int_0^Tu(\gamma^k)
\{u'(\gamma^k)^*(\frac{z^\delta-u(\gamma^k)}{u(\gamma^k)})\}dt\|_{\Gamma_i}^2+\beta_k
\|\gamma-\gamma^k\|_{\Gamma_i}^2\nb\\
&&+\Big\{ \int_{0}^T\|z^\delta-u(\gamma^{k})\|_{\Gamma_a}^2dt
-A\|\frac{1}{A}\int_0^Tu(\gamma^k)
\{u'(\gamma^k)^*(\frac{z^\delta-u(\gamma^k)}{u(\gamma^k)})\}dt\|_{\Gamma_i}^2\Big\}.
\label{che2t}
\eqn
We note that the last term in \eqref{che2t} is a constant, so it does not affect the minimization.
Hence we will drop that term in the functional $J^s(\gamma,\gamma^{k})$ and consider the following
minimization:
\beqn
\min_{\gamma\in K}J^s(\gamma,\gamma^{k})
=\min_{\gamma\in K}A\|\gamma-\gamma^{k}-\f{1}{A}\int_0^Tu(\gamma^k)
\{u'(\gamma^k)^*(\frac{z^\delta-u(\gamma^k)}{u(\gamma^k)})\}dt\|_{\Gamma_i}^2
+\beta_k\|\gamma-\gamma^k\|_{\Gamma_i}^2.
\label{che3t}
\eqn
This is a simple quadratic minimization, and it is easy to find its exact minimizer explicitly:
\beqn
\underset{\gamma\in K}{\rm argmin}\,J^s(\gamma,\gamma^{k})
=\gamma^k+\frac{1}{A+\beta_k}\int_0^Tu(\gamma^k)
\{u'(\gamma^k)^*(\frac{z^\delta-u(\gamma^k)}{u(\gamma^k)})\}dt.
\label{che5t}
\eqn
This motivates us with the following reconstruction algorithm for the Robin coefficient
in (\ref{ch1t}), which is clearly much easier and computationally much less expensive than
solving the minimization (\ref{ch4t}) directly.

\begin{algorithm}\label{al:par}

 Choose a tolerance parameter $\epsilon>0$ and an initial value
$\gamma^{0}$, set $k:=0$.
\bn
\item Compute $\gamma^{k+1}$:
\beqn
\gamma^{k+1}&=&\gamma^k+\frac{1}{A+\beta_k}\int_0^Tu(\gamma^k)
\{u'(\gamma^k)^*(\frac{z^\delta-u(\gamma^k)}{u(\gamma^k)})\}dt.
\label{cc6t}
\eqn

\item If $\frac{\|\gamma^{k+1}-\gamma^{k}\|_{\Gamma_i}}{\|\gamma^{k}\|_{\Gamma_i}}
\leq\epsilon$, stop the iteration;
otherwise set $k:=k+1$, go to Step 1.

\en
\end{algorithm}

%

\section{Numerical experiments }\label{sec:nsource}
\setcounter{equation}{0}
\setcounter{figure}{0}
\setcounter{table}{0}
In this section, we shall apply Algorithms \ref{al:ell} and \ref{al:par} that were proposed in
the previous Subsections~\ref{sec:ddmrobin} and \ref{sec:ddmrobint} to identify the Robin coefficients
in the elliptic and parabolic systems (\ref{ch1}) and (\ref{ch1t}) respectively.

We choose the domain $\Omega=(0,1)\times(0,2)$ and triangulate it  into
$N\times M$ small squares of equal size and further divide each
square through its diagonal into two triangles. This results in a finite element triangulation
of domain $\Om$.
All the elliptic problems involved in Algorithms \ref{al:ell} are solved by the continuous linear finite
element method, while all the parabolic problems in Algorithm \ref{al:par} are solved
by the continuous linear finite
element method in space and the backward difference scheme in time.

The parameters involved in Algorithms \ref{al:ell} and \ref{al:par} are chosen as follows.
The initial guesses are set to be identically equal to some constants,
which as we see are rather poor initial guesses for all the test problems.
The noisy data $z^\delta$ is obtained by adding
some uniform random noise to the exact data, i.e.,
$z^\delta=u + \delta R\, u$, where $R$ is a uniform random function varying in the
range [-1,1].
The errors are the relative $L^2$-norm errors
$\|\gamma^{k}-\gamma^*\|/\|\gamma^*\|$, where $\gamma^*$ and $\gamma^{k}$ are the exact parameter and
its numerical reconstruction by  Algorithms \ref{al:ell} and \ref{al:par} at the $k$th iteration respectively.

\medskip
We start two numerical tests for the Robin coefficient reconstructions on the partial
boundary $\Gamma_i=\{(x,y); ~x=1,\,\,0\leq y\leq2\}$ in the elliptic system
(\ref{ch1}), where we take
$a(\x)=c(\x)=1$ in $\Om$, the ambient temperature $g=2+(\cos(\pi y)+1)\gamma(\x)$ on $\Gamma_i$,
the heat flux $h=0$ on $\Gamma_a$, the source strength $f=(\pi^2+1)\cos(\pi y)+x^2-2$ and the exact forward solution
$u=x^2+\cos(\pi y)$ in $\Om$. We set the noise level $\delta=2\%$, the mesh $N=16$ and $M=32$,
the tolerance parameter $\epsilon=2\times 10^{-3}$ and the constant $A=1$.

\begin{example}
We take  the exact Robin coefficient $\gamma=3-\sin(\frac{\pi}{2}y)$ and
the initial guess $\gamma^{(0)}=2$.
\label{ex5}
\end{example}
%

%

\begin{figure}
  \centering
  \begin{tabular}{cc}
    \includegraphics[width=7cm]{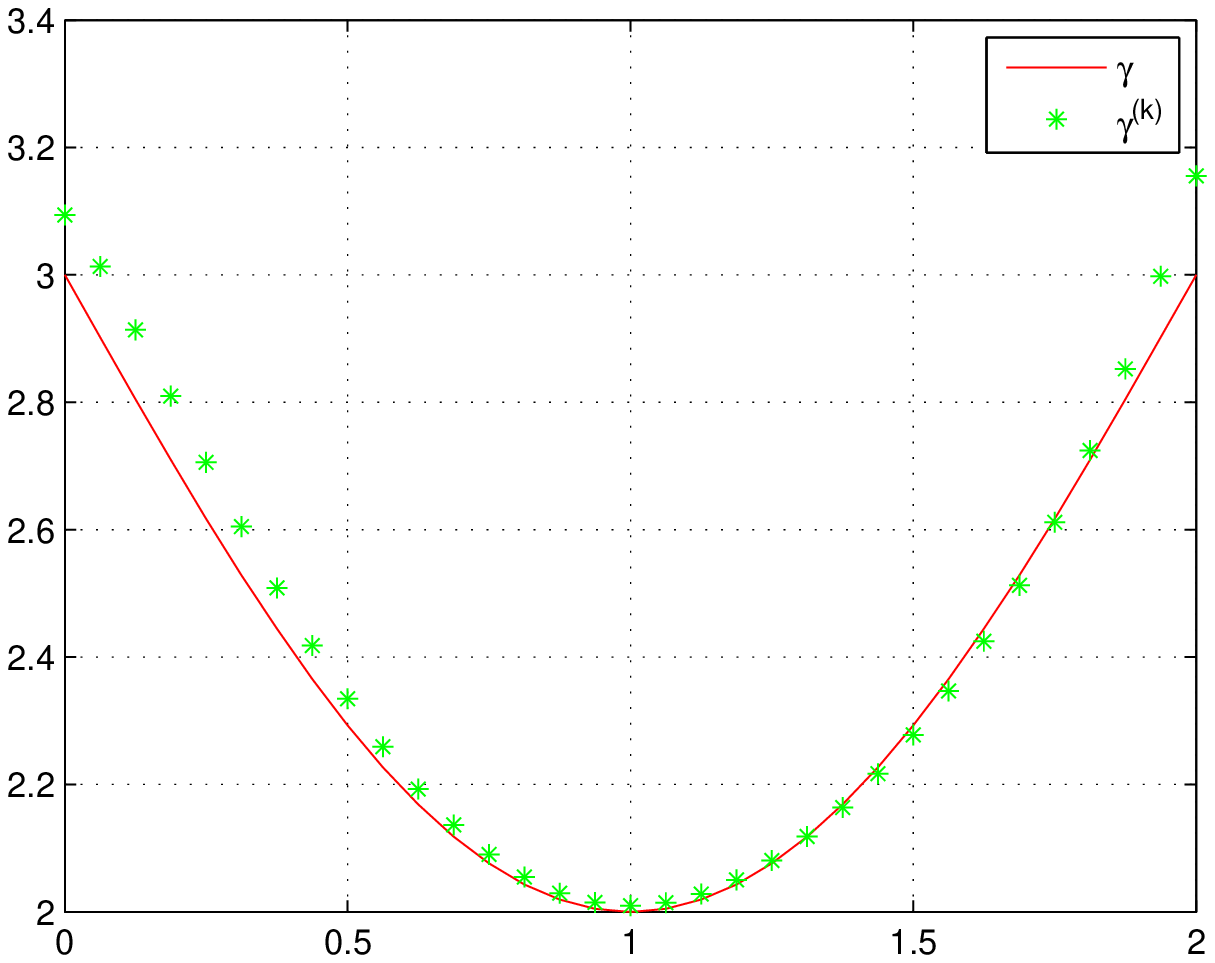} & \includegraphics[width=7cm]{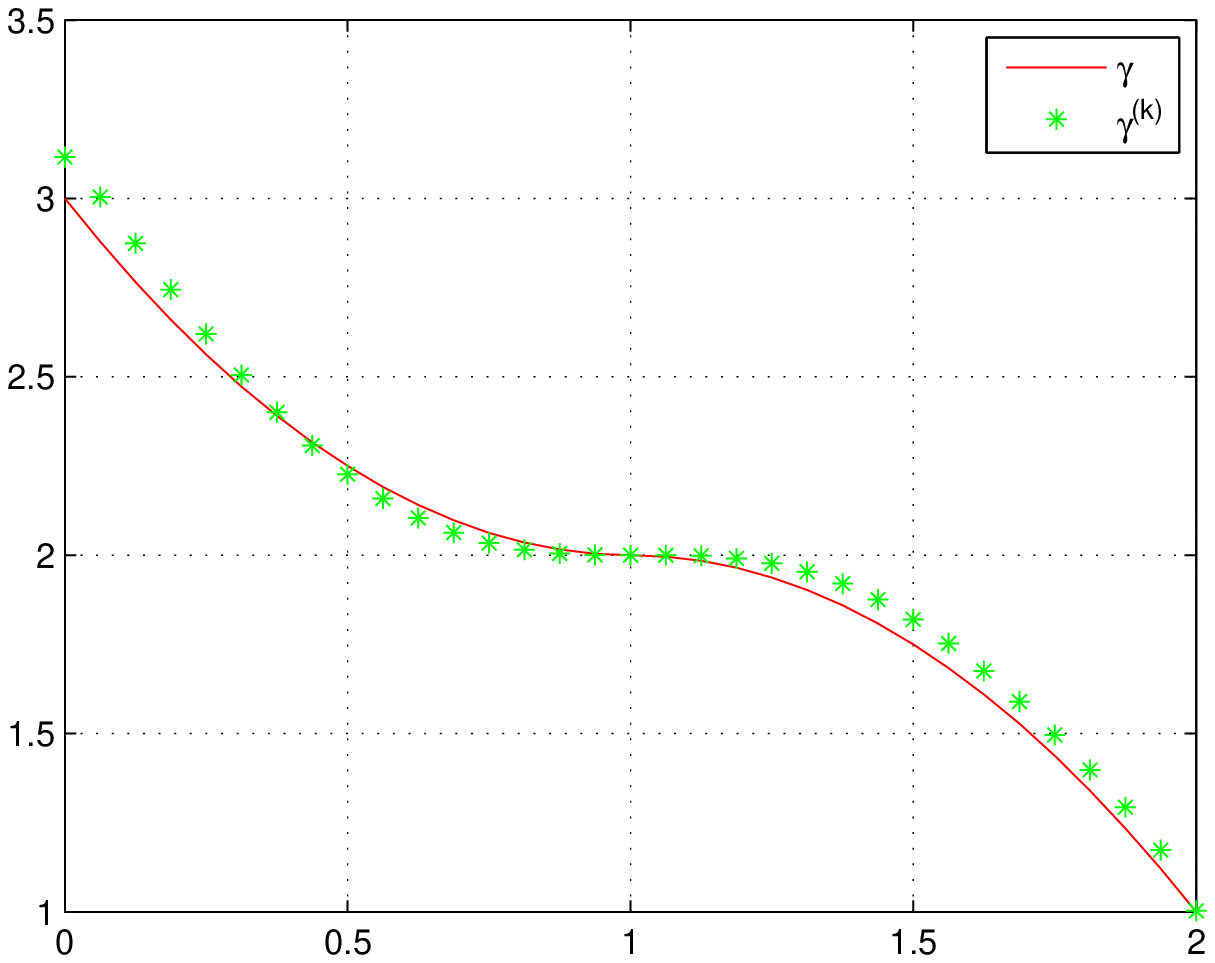} \\
  \end{tabular}
  \caption{Exact and reconstructed Robin coefficients for Example \ref{ex5} (left:
  $k=13$, the error $=0.0235$) and
  Example \ref{ex6} (right: $k=19$, the error $=0.0270$).}
  \label{aa2}
\end{figure}

\begin{example}
We take  the exact Robin coefficient $\gamma=(y-1)^2+2$ on $\{(x,y)\in \Gamma_i;\,
0\leq y\leq 1\}$ and $\gamma=-(y-1)^2+2$ on $\{(x,y)\in \Gamma_i;\,
1\leq y\leq 2\}$ and the initial guess $\gamma^{(0)}=2$.
\label{ex6}
\end{example}

Figure \ref{aa2} (left) and Figure \ref{aa2} (right)  give respectively the exact and reconstructed Robin
coefficients for Examples \ref{ex5} and \ref{ex6}.
We  see from the figure that the numerical reconstructed Robin coefficients, with a $2\%$ noise in the data
and very bad constant initial guesses, appear to be quite satisfactory, in view of
the severe ill-posedness of the inverse Robin problem. We can also clearly see
that Algorithms \ref{al:ell} converges quite fast with less than 20 iterations.

\medskip
Next, we demonstrate two numerical examples of reconstructing the Robin coefficient
$\gamma(\x)$ on the partial
boundary $\Gamma_i=\{(x,y); ~x=1,\,\,0\leq y\leq2\}$ in the parabolic system (\ref{ch1t}) with
$a(\x)=1$. We take the ambient temperature $g=(2+(\cos(\pi y)+1)\gamma(\x))t$
on $\Gamma_i\times [0,T]$, the heat
flux $h=0$ on $\Gamma_a\times [0,T]$, the source strength $f=\cos(\pi y)+x^2+(\pi^2\cos(\pi y)-2)t$ and
the exact forward solution $u=(x^2+\cos(\pi y))t$ in $\Om\times [0,T]$. We set the noise level $\delta=2\%$,
the mesh $N=16$ and $M=32$,
the tolerance parameter $\epsilon=5\times10^{-3}$,  the constant $A=1$ and the terminal time $T=2$.

\begin{example}
We take  the exact Robin coefficient $\gamma=-(y-1)^2+2$ and
 the initial guess $\gamma^{(0)}=2$.
\label{ex3}
\end{example}

\begin{figure}
  \centering
  \begin{tabular}{cc}
    \includegraphics[width=7cm]{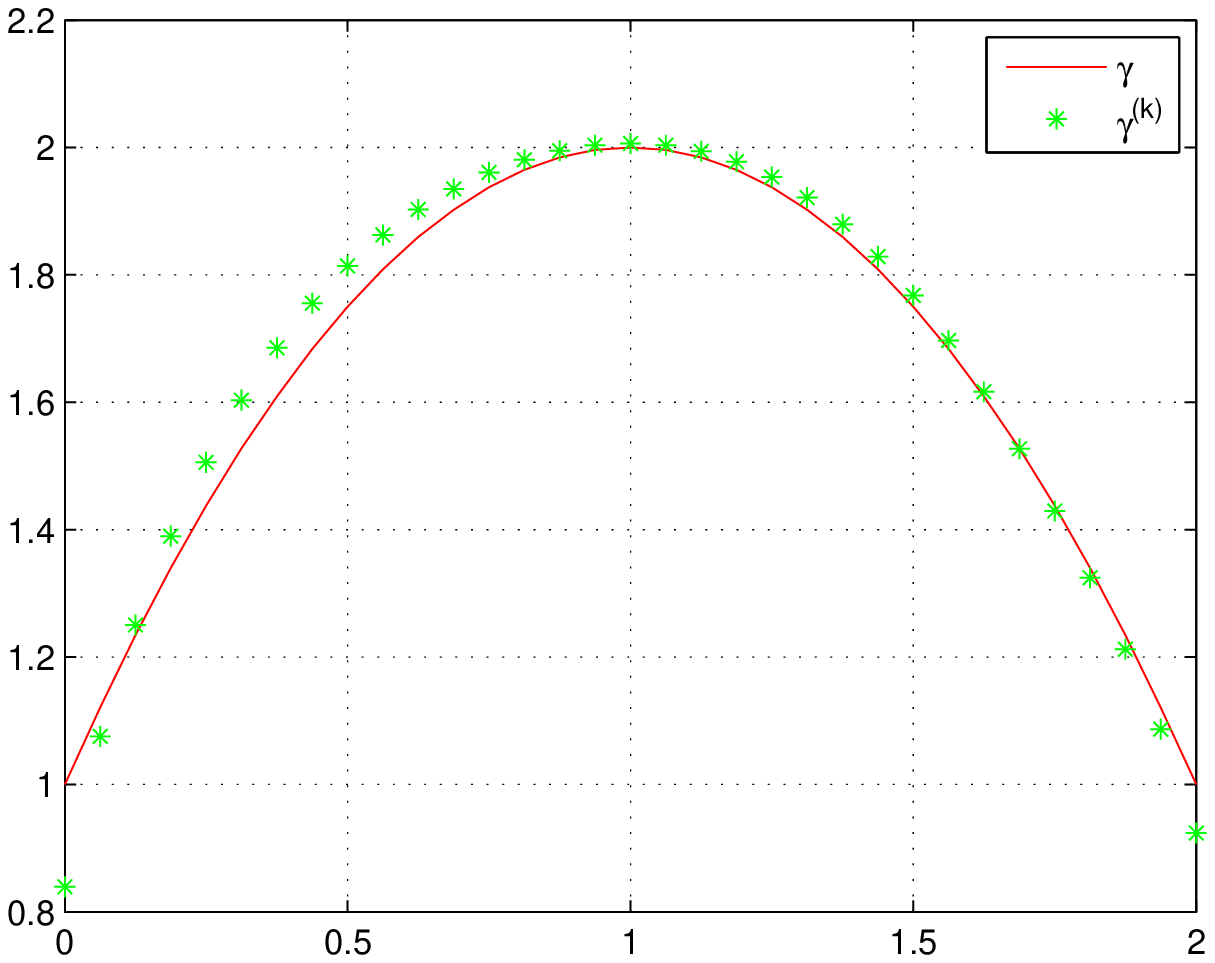} & \includegraphics[width=7cm]{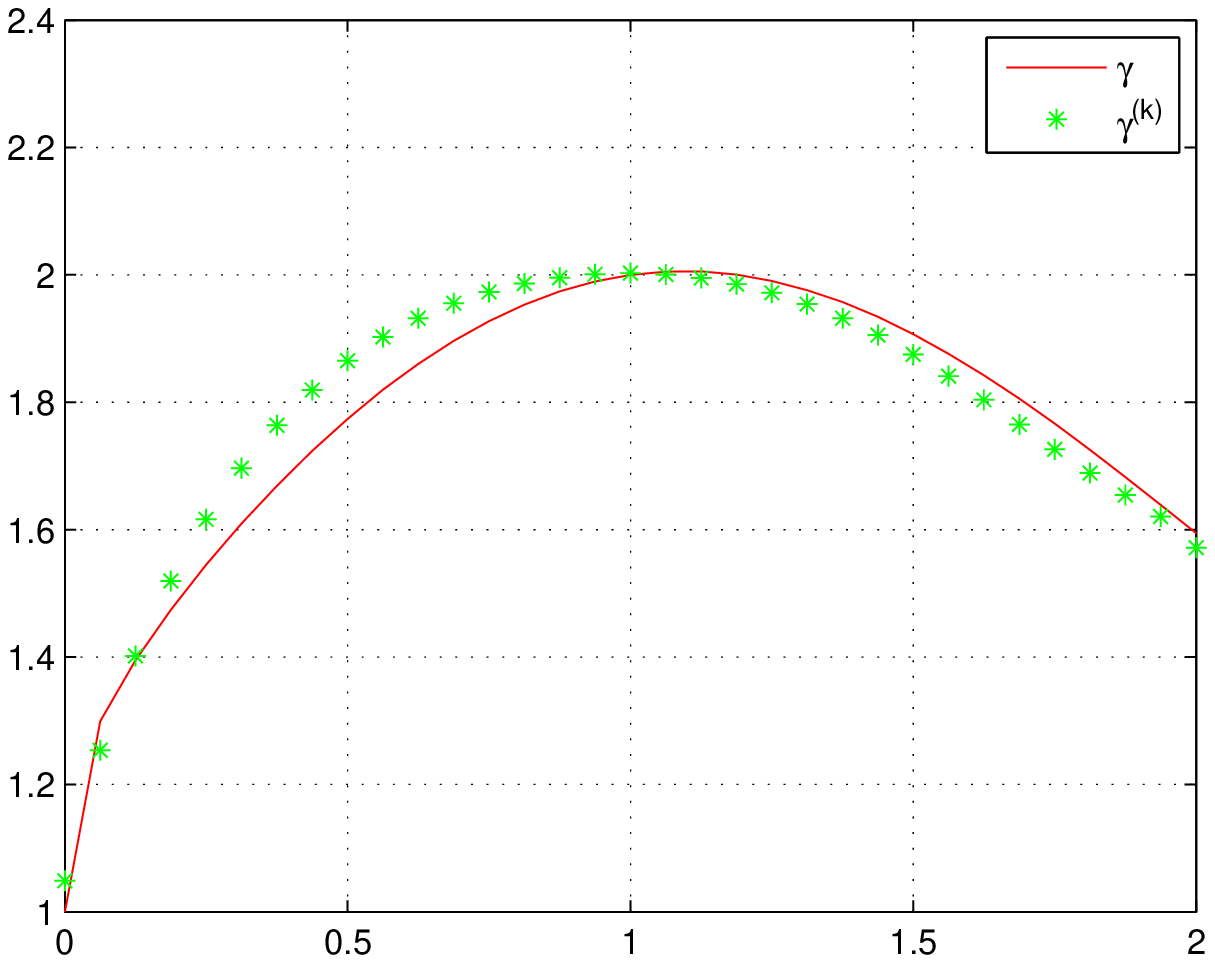} \\
  \end{tabular}
  \caption{Exact and reconstructed Robin coefficients for Example \ref{ex3} (left:
  $k=14$, the error $=0.0281$) and
  Example \ref{ex7} (right: $k=12$, the error $=0.0270$)}
\label{bb1}
\end{figure}

\begin{example}
We take  the exact Robin coefficient $\gamma=\frac{1}{2}(\sin(\frac{\pi}{2}y)+y^{\frac{1}{4}})+1$ and
 the initial guess $\gamma^{(0)}=2$.
\label{ex7}
\end{example}

Figure \ref{bb1} (left) and Figure \ref{bb1} (right)  give respectively the exact and reconstructed
Robin coefficients for Examples \ref{ex3} and \ref{ex7}.
%
%
%
%
%
We  see from  the figure \ref{bb1} that the
numerical reconstructed Robin coefficients, with a $2\%$ noise in the data
and very bad constant initial guesses, appear to be quite satisfactory, in view of
the severe ill-posedness of the inverse Robin problem. And we also note
that Algorithm \ref{al:par} converges quite fast with less than 20 iterations.

\section{Concluding remarks }\label{sec:conclu}
We have justified in this work the uniqueness of the elliptic and parabolic
Robin inverse problems.
Then the Levenberg-Marquardt iterative method is formulated to solve the nonlinear
Tikhonov regularized optimizations, which transform the original highly nonlinear and nonconvex
minimizations into convex minimizations.
We have established the quadratic convergence and the quadratic rate of convergence
for the L-M iterations
for the highly ill-posed nonlinear elliptic and parabolic Robin inverse problems.
This appears to be the first time in literature to achieve
the quadratic convergence and the quadratic rate of convergence for the L-M iterations rigorously for
a highly nonlinear and ill-posed inverse problem, in combination with a simple and easily implementable
choice rule of regularization parameters.
The surrogate functional techniques have been applied
to solve the convex minimizations at each L-M iteration,
which lead to explicit expressions of the minimizers for
both the elliptic and parabolic cases, resulting in two computationally very efficient
solvers for the highly ill-posed nonlinear inverse problems.
Numerical experiments have demonstrated the computational efficiency of the methods and
their robustness against the noise in the observation data.
%

\end{document}